\documentclass[10pt,a4paper]{article}
\usepackage{amsmath,amsthm,amssymb,bm,xspace,mydefs,diagrams,enumerate,color}
\renewcommand{\forceindent}{{}}
\usepackage[numbers]{natbib}
\begin{document}

\title{Analysis of the edge finite element approximation of the
  Maxwell equations with low regularity solutions} \author{Alexandre
  Ern\footnotemark[1] \and Jean-Luc Guermond\footnotemark[2]}
\footnotetext[1]{Universit\'e Paris-Est, CERMICS (ENPC), 77455
  Marne-la-Vall\'ee cedex 2, France and INRIA Paris, 
  75589 Paris, France.}  \footnotetext[2]{Department of Mathematics,
  Texas A\&M University 3368 TAMU, College Station, TX 77843, USA}

\maketitle

\begin{abstract}
We derive $\Hrt$-error estimates and improved $\bL^2$-error
estimates for the Maxwell equations approximated using edge finite elements. 
These estimates only invoke the
expected regularity pickup of the exact solution in the scale of the
Sobolev spaces, which is typically lower than $\frac12$ and can be
arbitrarily close to $0$ when the material properties are
heterogeneous. The key tools for the analysis are commuting
quasi-interpolation operators in $\Hrt$- and $\Hdv$-conforming
finite element spaces and, most crucially, newly-devised
quasi-interpolation operators delivering optimal estimates on the
decay rate of the best-approximation error for functions with
Sobolev smoothness index arbitrarily close to $0$. The
proposed analysis entirely bypasses the technique known in the
literature as the discrete compactness argument.
\end{abstract}


% \fi

%

\section{Introduction} 
The objective of this paper is to review some recent results
concerning the approximation of the Maxwell equations using edge
finite elements. One important difficulty is the modest regularity
pickup of the exact solution in the scale of the Sobolev spaces which
is typically lower than $\frac12$ and can be arbitrarily close to $0$
when the material properties are heterogeneous.  We show that the
difficulties induced by the lack of stability of the canonical
interpolation operators in $\Hrt$- and $\Hdv$-conforming finite
element spaces can be overcome by invoking recent results on commuting
quasi-interpolation operators and newly devised quasi-interpolation
operators that deliver optimal estimates on the decay rate of the
best-approximation error in those spaces.  In addition to a
curl-preserving lifting operator introduced by
\citet[p.~249-250]{Monk_1992}, the commuting quasi-interpolation
operators are central to establish a discrete counterpart of the
Poincar\'e--Steklov inequality (bounding the $\bL^2$-norm of a
divergence-free field by the $\bL^2$-norm of its curl), as already
shown in the pioneering work of~\cmag{\citet[\S9.1]{ArnFW:06}} on
Finite Element Exterior Calculus. It is therefore possible to bypass
entirely the technique known in the literature as the discrete
compactness argument
(\citet{Kikuchi_1989,Monk_Dem_2001,Caorsi_Fer_Raf_2000}).  The novelty
here is the use of quasi-interpolation operators devised by the
authors in~\cite{ErnGuermond:15} that give optimal decay rates of the
approximation error in fractional Sobolev spaces with a smoothness
index that can be arbitrarily small. This allows us to establish
optimal $\Hrt$-norm and $\bL^2$-norm error estimates that do not
invoke additional regularity assumptions on the exact solution other
than those resulting from the model problem at hand.  Optimality is
understood here in the sense of the decay rates with respect to the
mesh-size; the constants in the error estimates can depend on the
heterogeneity ratio of the material properties. Note that all
  the above quasi-interpolation operators are available with or
  without prescription of essential boundary conditions.

The paper is organized as follows. Notation and technical results are
given in \S\ref{Sec:Preliminaries}.  The main results from this
section are Theorem~\ref{Th:Quasi_interpolation_commutes}, which
states the existence of optimal commuting quasi-interpolation
operators, and Theorem~\ref{Th:glob_best_app} and
Theorem~\ref{Th:best_approx_Dir}, which give decay estimates of the
best approximation in fractional Sobolev norms.
\S\ref{Sec:Maxwell:model_problem} is concerned with standard facts
about the Maxwell equations. In particular, we state our main
assumptions on the model problem and briefly recall standard
approximation results for the Maxwell equations that solely rely on a
coercivity argument.  The new results announced above are collected in
\S\ref{Sec:coercivity:revisited} and in \S\ref{Sec:L2}. After
establishing the discrete Poincar\'e--Steklov inequality in
Theorem~\ref{Thm:Poincare}, our main results are
Theorem~\ref{Th:Hrot_convergence_MHD_Maxwell_robust} for the
$\Hrt$-error estimate and Theorem~\ref{Th:improved_L2_estimate} for
the improved $\bL^2$-error estimate. Both results do not invoke
regularity assumptions on the exact solution other than
those resulting from the model problem at hand.

\section{Preliminaries}
\label{Sec:Preliminaries}
We recall in this section some notions of functional analysis and
approximation using finite elements that will be invoked in the paper.
The space dimension is $3$ in the entire paper ($d=3$) and $\Dom$ is
an open, bounded, and connected Lipschitz subset in $\Real^3$.

\subsection{Functional spaces}
We are going to make use of the standard $L^2$-based Sobolev spaces
$H^m(\Dom)$, $m\in \Natural$.  The vector-valued counterpart of
$H^m(\Dom)$ is denoted $\bH^m(\Dom)$. We additionally introduce the
vector-valued spaces
\begin{align}
\Hrot &:= \{\bb\in \Ldeuxd\st \ROT \bb\in \Ldeuxd\}, \\
\Hdiv &:= \{\bb\in \Ldeuxd\st \DIV \bb\in \Ldeuxd\}.
\end{align}
 To be dimensionally
coherent, we equip theses Hilbert spaces with the norms
\begin{align}
\|b\|_{\Hun}    &:= (\|b\|_{\Ldeux}^2 +\ell_D^2\|\GRAD b\|_{\Ldeux}^2)^{\frac12},\\
\|\bb\|_{\Hrot} &:= (\|\bb\|_{\Ldeux}^2 +\ell_D^2\|\ROT\bb\|_{\Ldeux}^2)^{\frac12},\\
\|\bb\|_{\Hdiv} &:= (\|\bb\|_{\Ldeux}^2 +\ell_D^2\|\DIV\bb\|_{\Ldeux}^2)^{\frac12},
\end{align}
where $\ell_D$ is some characteristic dimension of $\Dom$, say the
diameter of $\Dom$ for instance.  In this paper we are also going to
use fractional Sobolev norms with smoothness index $s\in (0,1)$, defined as
follows:
\begin{equation}
 \|\bb\|_{\bH^s(\Dom)}:=(\|\bb\|_{\Ldeuxd}^2 +
\ell_\Dom^{2s} |\bb|_{\bH^s(\Dom)}^2)^{\frac12},
\end{equation} 
where $|\cdot|_{\bH^s(\Dom)}$ is the Sobolev--Slobodeckij semi-norm
applied componentwise. Similarly, for any $s>0$,
$s\in\Real\setminus\polN$, and $p\in[1,\infty)$, the norm of the
Sobolev space $W^{s,p}(\Dom)$ is defined by
$\|v\|_{W^{s,p}(\Dom)} := (\|v\|_{W^{m,p}(\Dom)}^p + \ell_D^{sp}
\sum_{|\alpha|=m}|\partial^\alpha
v|_{W^{\sigma,p}(\Dom)}^p)^{\frac1p}$
with
$\|v\|_{W^{m,p}(\Dom)} := (\sum_{|\alpha|\le m}
\ell_\Dom^{|\alpha|p}\|\partial^\alpha v \|_{L^p(\Dom)}^p)^{\frac1p}$
where $m:=\lfloor s\rfloor\in \polN$, $\sigma:=m-s\in(0,1)$.

\subsection{Traces} \label{sec:traces}
In order to make sense of the boundary conditions, we introduce trace operators.
Let $\gamma\upg:\Hun \to H^{\frac12}(\front)$ be the (full) trace
operator. It is known that $\gamma\upg$ is surjective.  Let
$\langle\SCAL,\SCAL\rangle_\front$ denote the duality pairing between
$\bH^{-\frac12}(\front):= (\bH^{\frac12}(\front))'$ and $\bH^{\frac12}(\front)$. 
We define the tangential trace operator $\gamma\upc:\Hrot\to
\bH^{-\frac12}(\front)$ as follows:
\begin{align}
\langle \gamma\upc(\bv), \bl\rangle_\front
  := \int_\Dom \bv\SCAL \ROT \bw(\bl)\dif x - \int_\Dom (\ROT\bv)\SCAL \bw(\bl)\dif x 
\label{eq:int_by_parts_rot_rot},
\end{align}
for all $\bv\in \Hrot$, all $\bl\in \bH^{\frac12}(\front)$ and all
$\bw(\bl)\in \bH^1(\Dom)$ such that $\gamma\upg(\bw(\bl))=\bl$.  One
readily verifies that the definition~\eqref{eq:int_by_parts_rot_rot}
is independent of the choice of $\bw(\bl)$,
% Indeed, let $\bw_1,\bw_2\in
% \bW^{1,p'}(\Dom)$ be such that $\gamma\upg(\bw_1)=\gamma\upg(\bw_2)=\bl$,
% \ie $\bw_1-\bw_2\in \bW^{1,p'}_0(\Dom)$.  Let $(\bphi_n)_{n\in
%   \Natural}$ be a sequence in $\bC_0^\infty(\Dom)$ converging to
% $\bw_1-\bw_2$ in $\bH^1_0(\Dom)$. Then, $0=\int_\Dom \bv\SCAL
% \ROT\bphi_n\dif x-\int_\Dom \bphi_n\SCAL \ROT \bv\dif x $.  Passing to
% the limit $n\to\infty$ yields $0=\int_\Dom \bv\SCAL
% \ROT(\bw_1-\bw_2)\dif x -\int_\Dom (\bw_1-\bw_2)\SCAL \ROT\bv\dif x $;
% hence, $\langle \gamma\upc(\bv), \gamma\upg(\bw_1)\rangle_\front
% =\langle \gamma\upc(\bv), \gamma\upg(\bw_2)\rangle_\front$,
% which establishes the claim.
that
$\gamma\upc(\bv)=\bv_{|\front}\CROSS\bn$ when $\bv$ is
smooth, and that the map $\gamma\upc$ is bounded.

We define similarly the normal trace map $\gamma\upd:\Hdiv \rightarrow
H^{-\frac12}(\front)$ by
\begin{align}
\langle \gamma\upd(\bv), l\rangle_\front
  := \int_\Dom \bv\SCAL \GRAD q(l)\dif x + \int_\Dom (\DIV \bv)q(l)\dif x, 
\label{eq:int_by_parts_div_grad} 
\end{align}
for all $\bv\in \Hdiv$, all $l\in H^{\frac12}(\front)$, and all
$q(l)\in H^1(\Dom)$ such that $\gamma\upg(q(l))=l$. Here
$\langle\SCAL,\SCAL\rangle_\front$ denotes the duality pairing between
$H^{-\frac12}(\front)$ and $H^{\frac12}(\front)$.  
One can verify that the
definition~\eqref{eq:int_by_parts_div_grad} is independent of the
choice of $q(l)$, that $\gamma\upd(\bv) = \bv_{|\front}\SCAL \bn$ when
$\bv$ is smooth, and that the map $\gamma\upd$ is bounded.

\begin{Th}[Kernel of trace operators]
  Let $\Hunz:=\overline{C^\infty_0(\Dom)}^{\Hun}$,
  $\Hrotz:= \overline{\bC^\infty_0(\Dom)}^{\Hrot}$, and
  $\Hdivz:= \overline{\bC^\infty_0(\Dom)}^{\Hdiv}$. Then, we have
\begin{subequations}\begin{align}
\Hunz =\KER(\gamma\upg), \\
\Hrotz = \KER(\gamma\upc), \\
\Hdivz = \KER(\gamma\upd).
\end{align}\end{subequations}
\end{Th}
\bproof 
The first identity is well-known, see \eg \citet[p.~315]{Brezis:11}.
The second one and the third one have been established in \cite[Thm.~4.9]{Ern_Guermond_CMAM_2016} in Lipschitz domains.
\eproof

\subsection{Generic finite element setting}
\label{Sec:Generic_FE_setting}
Let $\famTh$ be a shape-regular sequence of affine meshes.  To avoid
technical questions regarding hanging nodes, we also assume that the
meshes cover $\Dom$ exactly and that they are matching, \ie for all
cells $K,K'\in\calT_h$ such that $K\ne K'$ and $K\cap K'\ne\emptyset$,
the set $K\cap K'$ is a common vertex, edge, or face of both $K$ and
$K'$ (with obvious extensions in higher space dimensions).  Given a
mesh $\calT_h$, the elements in $K\in \calT_h$ are closed sets in
$\Real^3$ by convention, and they are all assumed to be constructed from a single reference cell $\wK$ through
affine, bijective, geometric transformations $\trans_K:\wK\to K$. 

The set of the mesh faces is denoted $\calF_h$ and is partitioned into
the subset of the interfaces denoted $\calFhi$ and the subset of the
boundary faces denoted $\calFhb$. Each interface $F$ is oriented by
choosing one unit vector $\bn_F$. The boundary faces are oriented by
using the outward normal vector. Given an interface $F\in \calFhi$, we
denote by $K_l$ and $K_r$ the two cells such that $F=K_l\cap K_r$ and
$\bn_F$ points from $K_l$ to $K_r$. This convention allows us to
define the notion of jump across $F$ as follows
\begin{equation}
\label{Eq:def_jump_scal}
\jump{v}_F(\bx) = v_{|K_l}(\bx) - v_{|K_r}(\bx) \qquad \textup{\ae $\bx$ in $F$}.
\end{equation}

We consider three types of reference elements in the sense of Ciarlet
as follows: $(\wK,\wP\upg,\wSigma\upg)$, $(\wK,\wbP\upc,\wSigma\upc)$
and $(\wK,\wbP\upd,\wSigma\upd)$. We think of
$(\wK,\wP\upg,\wSigma\upg)$ as a scalar-valued finite element with
some degrees of freedom that require point evaluations, for instance
$(\wK,\wP\upg,\wSigma\upg)$ could be a Lagrange element. The finite
element $(\wK,\wbP\upc,\wSigma\upc)$ is vector-valued with some
degrees of freedom that require to evaluate integrals over
edges. Typically, $(\wK,\wbP\upc,\wSigma\upc)$ is a N\'ed\'elec-type
or edge element.  Likewise, the finite element
$(\wK,\wbP\upd,\wSigma\upd)$ is vector-valued with some of degrees of
freedom that require evaluation of integrals over faces. Typically,
$(\wK,\wbP\upd,\wSigma\upd)$ is a Raviart-Thomas-type
element. The reader is referred to \citet{Hiptmair:99} for an
overview of a canonical construction of the above finite elements.

At this point we do not need to know the exact structure of the above
elements, but we are going to assume that they satisfy some commuting
properties. More precisely, let $s>\frac32$ and let
us consider the following functional spaces:
\begin{subequations}\begin{align}
\check V\upg(\wK) &=\{f\in H^{s}(\wK)\st \GRAD f\in \bH^{s-\frac12}(\wK)\}, \label{def:ck_Vg}\\
\check\bV\upc(\wK) &=\{\bg\in \bH^{s-\frac12}(\wK)\st \ROT \bg\in \bH^{s-1}(\wK)\}, \label{def:ck_Vc}\\
\check\bV\upd(\wK) &=\{\bg\in \bH^{s-1}(\wK)\st \DIV \bg\in L^1(\wK)\}.\label{def:ck_Vd}
\end{align}\end{subequations}
Let $\inter_{\wK}\upg, \inter_{\wK}\upc, \inter_{\wK}\upd$ be the
canonical interpolation operators associated with the above reference
elements.  Let $k\in\Natural$ and let $\polP_{k}(\Real^{3};\Real)$ be
the vector space composed of the trivariate polynomials of degree at
most $k$. We set $\wP\upb:=\polP_{k}(\Real^3;\Real)$ and let
$\inter_\wK\upb$ be the $L^2$-projection onto $\wP\upb$.  We now
state a key structural property that must be satisfied by the above Ciarlet
triples by assuming that the following diagram commutes:
\begin{equation}
\begin{diagram}[height=1.7\baselineskip,width=1.35cm]
    \check V\upg(\wK)  & \rTo^{\GRAD} & \check\bV\upc(\wK)  & \rTo^{\ROT} & \check\bV\upd(\wK)  & \rTo^{\DIV} & L^1(\wK)  \\
    \dTo_{\inter_\wK\upg} & & \dTo_{\inter_\wK\upc} & & \dTo_{\inter_\wK\upd}& & \dTo_{\inter_\wK\upb} \\
    \wP\upg & \rTo^{\GRAD} & \wbP\upc & \rTo^{\ROT} &
    \wbP\upd & \rTo^{\DIV} & \wP\upb 
\end{diagram}\label{Diag:wPupg_to_wPupb}
\end{equation}

In order to construct conforming approximation spaces based on $\famTh$ using the above
reference elements, we introduce the following linear maps:
\begin{subequations}\begin{align}
\mapKg(v) &=  v\circ\trans_K, \label{Eq:Def_mapg} \\
\mapKc(\bv) &= \Jac_K\tr (\bv\circ\trans_K), \label{Eq:def_Piola_rot}\\
\mapKd(\bv) &=  \det(\Jac_K)\,\Jac_K^{-1}(\bv\circ\trans_K),\label{Eq:Def_Piola} \\
\mapKb(v) &=\det(\Jac_K)(v\circ\trans_K),
\end{align}\end{subequations}
where $\mapKg$ is the pullback by $\trans_K$, and
$\mapKc$ and $\mapKd$ are the contravariant and covariant
Piola transformations, respectively. With these definitions in hand, we set
\begin{subequations} \label{eq:def_P} \begin{align}
P\upg(\calT_h)    &:=
\{v_h {\in} L^1(\Dom) \st \mapKg(v_{h|K}) {\,\in\,} \wP\upg, \ \forall K{\in\,} \calT_h, \
\jump{v_h}_F\upg=0,\ \forall F{\in}\calFhi\}, \hspace{-2pt} \label{def_Pupg}\\
\bP\upc(\calT_h)    &:= 
\{\bv_h {\in} \bL^1(\Dom) \st \mapKc(\bv_{h|K}) {\,\in\,} \wbP\upc, \ \forall K{\in\,} \calT_h, \
\jump{\bv_h}_F\upc=\bzero,\ \forall F{\in}\calFhi\}, \hspace{-2pt}\\
\bP\upd(\calT_h)    &:= 
\{\bv_h {\in} \bL^1(\Dom) \st \mapKd(\bv_{h|K}) {\,\in\,} \wbP\upd, \ \forall K{\in\,} \calT_h, \
\jump{\bv_h}_F\upd=0,\ \forall F{\in}\calFhi\},\hspace{-2pt}\\
P\upb(\calT_h)    &:= 
\{v_h\in L^1(\Dom) \st \mapKb(v_{h|K}) {\,\in\,} \wP\upb, \ \forall K{\in\,} \calT_h\},\label{def_Pupb}
\end{align} \end{subequations}
where $\jump{v_h}_F\upg:=\jump{v_h}_F$,
$\jump{\bv_h}_F\upc:=\jump{\bv_h}_F\CROSS \bn_F$, and
$\jump{\bv_h}_F\upd:=\jump{\bv_h}_F\SCAL\bn_F$.
Finally, to be able to account for boundary conditions, we define
\begin{subequations} \label{eq:def_P_0} \begin{align}
P\upg_0(\calT_h) &:= P\upg(\calT_h)\cap\Hunz,\label{def_Pupg_0}\\
\bP\upc_0(\calT_h) &:= \bP\upc(\calT_h)\cap\Hrotz,\\
\bP\upd_0(\calT_h) &:= \bP\upd(\calT_h)\cap\Hdivz. \label{def_Pupd_0}
\end{align} \end{subequations}

\subsection{Best approximation and commuting quasi-interpolation}
\label{Sec:Best_approximation_result}
Until recently, the stability and the convergence analysis of finite
element techniques for the approximation of the Maxwell equations was
made difficult by the absence of optimal approximation results. The
root of the difficulty was that the equivalent of the
Cl\'ement/Scott--Zhang quasi-interpolation operator was not available
for $\Hrot$-conforming and $\Hdiv$-conforming elements. Moreover, no
clear best-approximation estimate in fractional Sobolev norms was
known. We summarize in this section some of the most recent results in
this direction.

The bases for the construction of stable, commuting, and
quasi-interpolation projectors have been laid out in
\citet{Schoberl_2001,Schoberl_2008} and~\citet{Christ:07}, where
stability and commutation are achieved by composing the canonical
finite element interpolation operators with some mollification
technique. Then, following~\citet{Schoberl:05}, the projection
property over finite element spaces is obtained by composing these
operators with the inverse of their restriction to the said spaces.
An important extension of this construction allowing the possibility
of using shape-regular mesh sequences and boundary conditions has been
achieved by~\citet{Christiansen_Winther_mathcomp_2008} \cmag{(see also
\citet[\S5.4]{ArnFW:06} where this work was
prefigured)}. Further variants of this construction have lately been
proposed. For instance in~\citet{Christiansen:15}, the
quasi-interpolation projector has the additional property of
preserving polynomials locally, up to a certain degree, and
in~\citet{FalWi:12}, it is defined locally.  The results
of~\cite{Christiansen_Winther_mathcomp_2008} have been revisited in
\cite{Ern_Guermond_CMAM_2016} by invoking shrinking-based
mollification operators which do not require any extension outside the
domain.
% The results of \cite{Ern_Guermond_CMAM_2016} are
% stated in the language of numerical analysis, which makes them
% accessible to an audience who may not be familiar with the Finite
% Element Exterior Calculus language.

In order to stay general, we introduce an integer $q$ with the
convention that $q=1$ when we work with scalar-valued functions and
$q=3$ when we work with vector-valued functions.  For instance, we denote 
by $\polP_{k}(\Real^3;\Real^q)$ the vector space composed of the 3-variate polynomials
with values in $\Real^q$. The
quasi-interpolation results mentioned above can be summarized as
follows:
\begin{Th}[Stable, commuting projection] \label{Th:Quasi_interpolation_commutes} Let $P(\calT_h)$ be one of the finite element spaces
  introduced in \eqref{eq:def_P}-\eqref{eq:def_P_0}. Then there
  exists a quasi-interpolation operator
  $\calJ_h:L^1(\Dom;\Real^q) \to P(\calT_h)$ with the following
  properties:
\begin{enumerate}
\item[\textup{(i)}] $P(\calT_h)$ is pointwise invariant under $\calJ_h$.
\item[\textup{(ii)}] Let $p\in [1,\infty]$. There is $c$, uniform w.r.t. $h$, such that
  $\|\calJ_h\|_{\calL(L^p;L^p)}\le c$ and
\begin{equation}
\|f - \calJ_h (f) \|_{L^p(\Dom;\Real^q)} \le c \inf_{f_h \in P(\calT_h)} \|f -f_h\|_{L^p(\Dom;\Real^q)},
\label{Lp_error_estimate_for_calJh}
\end{equation}
for all $f\in L^p(\Dom;\Real^q)$;
\item[\textup{(iii)}] \label{Th:Pih_item3} $\calJ_h$ commutes with the standard
  differential operators, \ie the following diagrams commute:
\end{enumerate}
\begin{equation}
\begin{diagram}[height=1.7\baselineskip,width=1.5cm]
\Hun & \rTo^{\GRAD} & \Hrot & \rTo^{\ROT} & \Hdiv & \rTo^{\DIV} &L^2(\Dom)\\
\dTo_{\calJ\upg_h} & & \dTo_{\calJ\upc_h} & & \dTo_{\calJ\upd_h}& & \dTo_{\calJ\upb_h} \\
P\upg(\calT_h) & \rTo^{\GRAD} & \bP\upc(\calT_h) & \rTo^{\ROT} &
\bP\upd(\calT_h) & \rTo^{\DIV} & P\upb(\calT_h)  
\end{diagram}\label{Diag:calJh_commutes}
\end{equation}
\begin{equation}
\begin{diagram}[height=1.7\baselineskip,width=1.5cm]
\Hunz& \rTo^{\GRAD} & \Hrotz & \rTo^{\ROT} & \Hdivz & \rTo^{\DIV} &L^2(\Dom)\\
\dTo_{\calJ_{h0}\upg} & & \dTo_{\calJ_{h0}\upc} & & \dTo_{\calJ_{h0}\upd}& & \dTo_{\calJ_{h}\upb} \\
P\upg_0(\calT_h) & \rTo^{\GRAD} & \bP\upc_0(\calT_h) & \rTo^{\ROT} & \bP_0\upd(\calT_h) & \rTo^{\DIV} & P\upb(\calT_h) 
\end{diagram}\label{Diag:calJh0_commutes}
\end{equation}
\end{Th}

\begin{proof}
See \eg \cite[Thm.~6.5]{Ern_Guermond_CMAM_2016}.
\end{proof}

For the estimate \eqref{Lp_error_estimate_for_calJh} to be useful we
need an estimate on the decay rate of best-approximation error. 
This question is answered by the
following two results:
\begin{Th}[Best approximation] \label{Th:glob_best_app} Let $P(\calT_h)$ be one of the finite element spaces
introduced in \eqref{eq:def_P}. Let $k$ be the
largest integer such that
$\polP_{k}(\Real^3;\Real^q) \subset \wP$.  There exists a uniform
constant $c$ such that
\begin{equation} \label{Eq:Cor:interpolation_RT_N}
\inf_{w_h\in P(\calT_h)} \|v -w_h \|_{L^{p}(\Dom;\Real^q)}
\le c\, h^{r} |v|_{W^{r,p}(\Dom;\Real^q)},
\end{equation}
for all $r\in [0,k+1]$, all $p\in [1,\infty)$ if $r\not\in\polN$ or all $p\in [1,\infty]$ if $r\in\polN$, and all $v\in W^{r,p}(\Dom;\Real^q)$.
\end{Th}
\begin{proof}
See \cite[Cor.~5.4]{ErnGuermond:15}.
\end{proof}

\begin{Th}[Best approximation with boundary conditions] 
  \label{Th:best_approx_Dir}
  Let $P_0(\calT_h)$ be one of the finite element spaces
  introduced in \eqref{eq:def_P_0} and let $\gamma$ be the trace operator
from~\textup{\S\ref{sec:traces}} associated with $P_0(\calT_h)$. Let $k$ be the largest integer
  such that $\polP_{k}(\Real^3;\Real^q) \subset \wP$.  There exists a
  uniform constant $c$, that depends on $|rp-1|$, such that
\begin{equation}
\!\,\inf_{w_h \in P_0(\calT_h)}\!\|v-w_h\|_{L^p(\Dom;\Real^q)} \! \le 
\begin{cases}
c h^r |v|_{W^{r,p}(\Dom;\Real^q)},& \text{$\forall v\in W^{r,p}_{0,\gamma}(\Dom;\Real^q)$ if $rp>1$},\\
c h^r \ell_\Dom^{-r} \|v\|_{W^{r,p}(\Dom;\Real^q)},& \text{$\forall v\in W^{r,p}(\Dom;\Real^q)$ if $rp<1$},
\end{cases}\hspace{-0.5cm}
\end{equation}
where
$W^{r,p}_{0,\gamma}(\Dom;\Real^q):=\{v\in
W^{r,p}(\Dom;\Real^q)\st\gamma(v)=0\}$,
for all $r\in [0,k+1]$, all $p\in [1,\infty)$ if $r\not\in\polN$ or
all $p\in [1,\infty]$ if $r\in\polN$.
\end{Th}
\begin{proof}
See \cite[Cor.~6.5]{ErnGuermond:15}.
\end{proof}

Localized versions of the above results and best-approximation error
estimates for higher-order Sobolev semi-norms can be found in
\cite{ErnGuermond:15}. These results are proved by constructing
quasi-interpolation operators in a unified way. This construction is
done in two steps and consists of composing an elementwise
projection onto the broken finite element space with a smoothing
operator based on the averaging of the degrees of freedom on the
broken space.
We also refer the reader to \citet{Ciarlet:13} for similar estimates
for the Scott--Zhang quasi-interpolation operator in the context of
scalar-valued finite elements and $rp>1$.

\bRem[Edge elements] To put the above results in perspective with the
literature,
let us observe that the canonical interpolation operator for edge
elements is only stable in $\bH^s(\Dom)$ for $s>1$. Using the
techniques in \citet{AmBDG:98}, one can show that this operator is
also stable in the space
$\{\bv \in \bH^s(\Dom)\tq \ROT\bv \in \bL^p(\Dom)\}$ with $s>\frac12$
and $p>2$, see, e.g., \citet{BofGa:06}. In contrast with these
results, Theorem~\ref{Th:glob_best_app} states that any function $\bv$
in $\bH^s(\Dom)$, with $s$ arbitrarily close to
zero, can be optimally approximated in $\bP\upc(\calT_h)$. 
Let us also observe that the best-approximation result from
Theorem~\ref{Th:best_approx_Dir} can be used with functions that are not
smooth enough to have a well-defined trace at the boundary and that
these functions are approximated by finite element functions that do satisfy a
boundary condition.  
\eRem

\section{Maxwell's equations}
\label{Sec:Maxwell:model_problem} In this section we recall standard
facts about the Maxwell equations that will be used later in the
paper. For an introduction to the subject, the reader is referred to
\citet[Chap.~1]{Bossavit_GB} or \citet[Chap.~1]{Monk_2003}.

\subsection{The model problem}
Maxwell's equations consist of a set of partial differential equations
giving a macroscopic description of electromagnetic phenomena. More
precisely, these equations describe how the electric field, $\bE$, the
magnetic field $\bH$, the electric displacement field, $\bD$, and the
magnetic induction (sometimes called magnetic flux density), $\bB$,
interact through the action of currents, $\bj$, and charges, $\rho$:
\begin{subequations}\label{Maxwell_PDEs}\begin{alignat}{2}
&\partial_t \bD -\ROT \bH = -\bj &\quad&\text{(Amp\`ere's law),}\label{eq:Ampere}\\
&\partial_t \bB + \ROT \bE = \bzero&\quad&\text{(Faraday's law of induction),}\label{eq:Faraday}\\
&\DIV \bD=\rho &\quad&\text{(Gauss' law for electricity),}\label{eq:Gauss_e} \\
&\DIV \bB =0 &\quad&\text{(Gauss' law for magnetism).}\label{eq:Gauss_m}
\end{alignat}\end{subequations}
Note that~\eqref{eq:Gauss_e}-\eqref{eq:Gauss_m} can be viewed as
constraints on the time-evolution
problem~\eqref{eq:Ampere}-\eqref{eq:Faraday}. If $(\DIV\bB)_{|t=0}=0$,
then taking the divergence of~\eqref{eq:Faraday} implies that
\eqref{eq:Gauss_m} holds at all times. Similarly, assuming that
$(\DIV\bD)_{|t=0}=\rho_{|t=0}$ and that the charge conservation
equation $\partial_t\rho + \DIV\bj=0$ holds at all times implies that
\eqref{eq:Gauss_e} holds at all times.

The system~\eqref{Maxwell_PDEs} is closed by relating the fields through constitutive
laws describing microscopic mechanisms of polarization and
magnetization:
\begin{equation}
\bD -\varepsilon_0 \bE = \bP,\qquad \bB = \mu_0(\bH +\bM),
\end{equation}
where $\varepsilon_0$ and $\mu_0$ are the electric permittivity and
the magnetic permeability of vacuum; the fields $\bP$ and $\bM$ are
the polarization and the magnetization, respectively.  These
quantities are the average representatives at macroscopic scale of
complicated microscopic interactions, \ie they need to be modeled and
measured. For instance, $\bP=\bzero$ and $\bM=\bzero$ in vacuum, and it is
common to use $\bP=\varepsilon_0 \varepsilon_r \bE$ and $\bM=\mu_r\bH$
to model isotropic homogeneous dielectric and magnetic materials,
where $\varepsilon_r$ is the electric susceptibility and $\mu_r$ is
the magnetic susceptibility.  The current and charge density, $\bj$
and $\rho$, are a priori given, but it is also possible to make these
quantities depend on the other fields through phenomenological
mechanisms.  For instance, it is possible to further decompose the current
into one component that depends on the material and another one that
is a source; the simplest model doing that is Ohm's law: $\bj =\bj_s
+\sigma \bE$; $\sigma$ is the electrical conductivity and $\bj_s$ is
an imposed current.

We now formulate Maxwell's equations in two different regimes: the
harmonic regime, leading to the Helmholtz problem, and the eddy current limit.  
We henceforth assume that 
\begin{equation} \label{eq:DeE_BmH}
\bD =\epsilon \bE \quad\text{and}\quad \bB=\mu
\bH,
\end{equation} 
where $\epsilon$ and $\mu$ may be space-dependent.

\subsection{The Helmholtz problem}\label{Sec:Mawell_harminic}
We first consider Maxwell's equations in the harmonic regime. Using
the convention $i^2=-1$, the time-dependence is assumed to be of the
type $\text{e}^{i\omega t}$ where $\omega$ is the angular
frequency. Letting $(\front\subD,\front\subN)$ be a partition of the
boundary $\front$ of $\Dom$, the time-harmonic version of
\eqref{eq:Ampere}-\eqref{eq:Faraday} is
\begin{subequations}\label{wave_Maxwell}
\begin{alignat}{2}
&i\omega \epsilon \bE + \sigma \bE - \ROT \bH = -\bj_s, &\quad& \text{in $\Dom$}, \\ 
&i\omega \mu \bH + \ROT \bE
  = \bzero, &\quad& \text{in $\Dom$}, \\ 
&\bH_{|\front\subD}{\times}\bn =
  \ba\subD,\quad \bE_{|\front\subN}{\times}\bn=\ba\subN, &\quad& \text{on $\front$}.
\end{alignat} 
\end{subequations}
The dependent variables are the electric field, $\bE$, and the
magnetic field, $\bH$. The data are the conductivity, $\sigma$, the
permittivity, $\epsilon$, the permeability, $\mu$, the current,
$\bj_s$, and the boundary data $\ba\subD$ and $\ba\subN$.  The
material coefficients $\epsilon$ and $\mu$ are allowed to be
complex-valued.  The system \eqref{wave_Maxwell} models, for instance,
a microwave oven; see \eg \cite[Chap.~9]{Bossavit_GB}. The conditions
$\bH_{|\front\subD}{\times}\bn=\bzero$ and
$\bE_{|\front\subN}{\times}\bn=\bzero$ are usually called perfect conductor
and perfect magnetic conductor boundary conditions, respectively.

Let us assume that the modulus of the magnetic permeability $\mu$ is
bounded away from zero uniformly in $\Dom$. It is then
possible to eliminate $\bH$ by using $\bH = i (\omega \mu)^{-1} \ROT
\bE$, and the resulting system takes the following form:
\begin{subequations} \label{eq:wave_Maxwell}\begin{alignat}{2}
&(-\omega^2 \epsilon + i\omega\sigma)\bE + \ROT(\mu^{-1} \ROT \bE) = -i\omega \bj_s,&\quad& \text{\ae in $\Dom$},\\
&(\ROT \bE)_{|\front\subD}\CROSS \bn = -i\omega \mu \ba\subD,\quad
 \bE_{|\front\subN}\CROSS \bn=\ba\subN,&\quad& \text{\ae on $\front$}.
\end{alignat}\end{subequations} 
If $\sigma=0$, this problem leads to two different situations
depending on whether $\omega$ is a resonance frequency of the domain
$\Dom$ or not.  If it is the case, the above problem is an eigenvalue
problem, otherwise it is a boundary-value problem.

\subsection{Eddy current problem}
\label{Sec:Maxwell_eddy_current}
When the time scale of interest, say $\tau$, is such that the ratio
$\epsilon/(\tau \sigma)$ is very small, \ie $\epsilon/(\tau \sigma)\ll
1$, it is legitimate to neglect the so-called displacement current
density (\ie Maxwell's correction, $\partial_t\bD$).  This situation
occurs in particular in systems with moving parts (either solid or
fluids) with a characteristic speed that is significantly slower than
the speed of light. Assuming again that $(\front\subD,\front\subN)$
forms a partition of the boundary $\front$ of $\Dom$, the resulting
system takes the following form:
\begin{subequations}\label{MHD_Maxwell}
\begin{alignat}{2}
&\sigma \bE - \ROT \bH = -\bj_s, &\quad& \text{in $\Dom$}, \\ 
&\partial_t (\mu \bH) + \ROT \bE
  = \bzero, &\quad& \text{in $\Dom$}, \\ 
&\bH_{|\front\subD}{\times}\bn =
  \ba\subD,\quad \bE_{|\front\subN}{\times}\bn=\ba\subN, &\quad& \text{on
  $\front$}.
\end{alignat}
\end{subequations}
The system \eqref{MHD_Maxwell} arises in magneto-hydrodynamics (MHD);
in this case, $\bj_s$ is further decomposed into $\bj_s = \bj_s' 
+\sigma \bu\CROSS \bB$  where $\bu$ is the velocity of the fluid 
occupying the domain $\Dom$, \ie the actual current is  
decomposed into $\bj=\bj_s' + \sigma (\bE+\bu\CROSS\bB)$.

Let us assume that the conductivity $\sigma$ is bounded uniformly from
below in the domain $\Dom$ by a positive constant $\sigma_{\min}$.
It is then possible to eliminate the electric field from
\eqref{MHD_Maxwell} by using $\bE=\sigma^{-1}(\ROT \bH - \bj_s)$.  The new
system to be solved is re-written as follows:
\begin{subequations}\begin{alignat}{2}
&\partial_t (\mu \bH) +
\ROT (\sigma^{-1}\ROT \bH) = \ROT (\sigma^{-1} \bj_s),&\quad& \text{in $\Dom$},\\
&\bH_{|\front\subD}\CROSS \bn=\ba\subD, 
\quad (\sigma^{-1}\ROT \bH)_{|\front\subN} \CROSS \bn
= \bc\subN,&\quad& \text{on $\front$},
\end{alignat}\end{subequations}
where $\bc\subN=\ba\subN + \sigma^{-1}\bj_{s|\front\subN}\CROSS \bn$.
At this point, it is possible to further simplify the problem by
assuming that either the time evolution is harmonic,
$\bH(\bx,t)=\bH_{\text{\rm sp}}(\bx) e^{i\omega t}$, or the time
derivative is approximated as follows $\partial_t \bH(\bx,t)\approx
(\dt)^{-1}(\bH(\bx,t) -\bH(\bx,t-\dt))$, where $\dt$ is the time step
of the time discretization. The above system then reduces to solving
the following problem:
\begin{subequations}\label{eq:MHD_Maxwell}\begin{alignat}{2}
&\tmu \bH +
\ROT (\sigma^{-1}\ROT \bH) = \bbf,&\quad& \text{in $\Dom$},\\
&\bH_{|\front\subD}\CROSS \bn=\ba\subD, 
\quad (\sigma^{-1}\ROT \bH)_{|\front\subN}\CROSS \bn = \bc\subN,&\quad& \text{on $\front$},
\end{alignat}\end{subequations}
after appropriately renaming  the dependent variable and the data, say
either $\widetilde\mu := i\omega \mu$ and $\bbf=\ROT
(\sigma^{-1} \bj_s) $, or $\widetilde\mu := \mu (\dt)^{-1}$ and
$\bbf:=\ROT (\sigma^{-1} \bj_s) + \widetilde\mu \bH(x,t-\dt)$, \etc; note
that $\DIV\bbf=0$ in both cases.

\subsection{Abstract problem}
\label{Sec:Weak_MHD_Maxwell} 
The Helmholtz problem and the eddy current problem have a very
similar structure.  For simplicity, we restrict the scope to Dirichlet
boundary conditions; the techniques presented below can be adapted to 
handle Neumann boundary conditions as well. After lifting the Dirichlet 
boundary condition
and making appropriate changes of notation, the above two problems
\eqref{eq:wave_Maxwell} and \eqref{eq:MHD_Maxwell} can be reformulated
in the following common form: Find $\bA$ such that
\begin{equation}
  \widetilde\mu \bA + \ROT(\kappa \ROT \bA) = \bbf,\qquad 
  \bA_{|\front}\CROSS \bn=\bzero.
%\quad (\kappa\ROT \bA)_{|\front\subN}\CROSS \bn=\bzero.
\label{Strong_MHD_Maxwell}
\end{equation} 
% We have assumed that the Neumann data is zero to avoid unnecessary
% technicalities.  
We henceforth
assume that $\bef\in \bL^2(\Dom)$ and that $\DIV\bef=0$. Note
that, by taking the divergence of the PDE
in~\eqref{Strong_MHD_Maxwell}, this property implies that
\begin{equation} \label{eq:divA}
\DIV(\widetilde\mu \bA) = 0.
\end{equation}
This additional condition on $\bA$ plays a key role in
\S\ref{Sec:coercivity:revisited} and in \S\ref{Sec:L2}. 
Concerning the material properties $\tmu$ and $\kappa$, we make three
assumptions:
(i) Boundedness: $\widetilde\mu, \kappa
\in L^\infty(\Dom;\polC)$ and we
set $\mu_\sharp =\esssup_{\bx\in \Dom}|\widetilde \mu(\bx)|$ and
$\kappa_\sharp = \esssup_{\bx\in \Dom}|\kappa(\bx)|$;
(ii) Rotated positivity: there are real numbers $\theta$, $\mu_\flat>0$, and
$\kappa_\flat>0$ so that
\begin{equation} \label{Hyp:mu_kappa_Weak_MHD_Maxwell}
\infess_{\bx\in\Dom}\Re(\text{e}^{i\theta}\widetilde\mu(\bx)) \ge \mu_\flat\quad
\text{and}
\quad \infess_{\bx\in\Dom}\Re(\text{e}^{i\theta}\kappa(\bx)) \ge \kappa_\flat.
\end{equation}
We define the heterogeneity ratios
$\mu_{\sharp/\flat}:=\frac{\mu_\sharp}{\mu_\flat}$
and $\kappa_{\sharp/\flat}:= \frac{\kappa_\sharp}{\kappa_\flat}$.
We also define the magnetic Reynolds number
$\gamma_{\tmu,\kappa}=\mu_\sharp\ell_\Dom^2\kappa_\sharp^{-1}$. (Obviously,
if the material is highly contrasted, several magnetic Reynolds numbers
can be defined.)
(iii) Piecewise smoothness: there is a partition of
$\Dom$ into $M$ disjoint Lipschitz subdomains
$\Dom_1,\cdots,\Dom_M$ so that $\tmu_{|\Dom_i} \in W^{1,\infty}(D_i)$
and $\kappa_{|\Dom_i} \in W^{1,\infty}(D_i)$ for all $1\le i\le M$.

\bRem[Assumption~\eqref{Hyp:mu_kappa_Weak_MHD_Maxwell}] This
assumption holds, \eg if the coefficient $\kappa$ is real and if
$\tmu=\rho_\mu e^{i\theta_\mu}$ with
$\essinf_{\bx\in \Dom}\rho_\mu(\bx)=\rho_\flat>0$ and
$\theta_\mu(\bx) \in [\theta_{\min},\theta_{\max}] \subset (-\pi,
\pi)$
\ae in $\Dom$ with $\delta:=\theta_{\max} -\theta_{\min}< \pi$; indeed,
one can take $\theta =-\frac{1}{2}(\theta_{\min}+\theta_{\max})
\frac{\pi}{2\pi-\delta}$,
$\mu_\flat=\min(\cos(\theta_{\min}+\theta),\cos(\theta_{\max}+\theta))
\rho_\flat$ and
$\kappa_\flat = \cos(\theta)
\essinf_{\bx\in \Dom} \kappa(\bx)$.
For instance, this is the case for the Helmholtz problem (where
$\widetilde\mu = -\omega^2 \epsilon + i\omega\sigma$ and
$\kappa = \mu^{-1}$) and for the eddy current problem (where
$\widetilde\mu := i\omega \mu$ or $\widetilde\mu := \mu (\dt)^{-1}$
and $\kappa = \sigma^{-1}$).  An important example when the condition
\eqref{Hyp:mu_kappa_Weak_MHD_Maxwell} does not hold is when the two
complex numbers $\tmu$ and $\kappa$ are collinear and point in
opposite directions; in this case, \eqref{Strong_MHD_Maxwell} is an
eigenvalue problem.  
\eRem
 
% \reply{There is nothing wrong with this calculation. However, there
% is one more condition to be verified, namely that
% $\theta\in (-\frac{\pi}{2},\frac{\pi}{2})$, so that
% $\Re(e^{i\theta}\kappa)$ also remains positive, and this condition
% is not likely to hold with the above choice. This is why I suggest
% to take
% \[
% \theta = -\mu \frac{\pi}{2\pi-\delta}, \qquad \mu=\frac{\theta_{\max}+\theta_{\min}}{2}, \qquad
% \delta = \theta_{\max} -\theta_{\min}.
% \]
% Then we have
% \begin{itemize}
% \item $\theta<\frac{\pi}{2}$ iff $-2\mu<2\pi-\delta$ iff $-\theta_{\min}<\pi$
% \item $\theta>-\frac{\pi}{2}$ iff $2\mu<2\pi-\delta$ iff $\theta_{\max}<\pi$
% \item $\theta_{\max}+\theta<\frac{\pi}{2}$ iff $\frac{\delta}{2}+\mu\frac{\pi-\delta}{2\pi-\delta}<\frac{\pi}{2}$ iff $2\mu<2\pi-\delta$ iff $\theta_{\max}<\pi$
% \item $\theta_{\min}+\theta>-\frac{\pi}{2}$ iff $-\frac{\delta}{2}+\mu\frac{\pi-\delta}{2\pi-\delta}>-\frac{\pi}{2}$ iff $-2\mu<2\pi-\delta$ iff $-\theta_{\min}<\pi$
% \end{itemize}
% Note that for the previous example where $\theta_{\min}=\theta_{\max}=\frac{\pi}{2}$, we obtain $\mu=\frac{\pi}{2}$, $\delta=0$ and $\theta=\frac{\pi}{4}$ as announced, whereas $\theta= - \frac{\theta_{\max}+\theta_{\min}}{2}=\frac{\pi}{2}$ is not legitimate.
% }

\subsection{Basic weak formulation and approximation}
\label{Sec:basic}
We recall in this section standard approximation results for
\eqref{Strong_MHD_Maxwell} that solely rely on a coercivity
argument. More subtle arguments are invoked
in~\S\ref{Sec:coercivity:revisited} and in~\S\ref{Sec:L2}; in
particular, we do not use here that $\DIV\bef=0$.

A weak formulation of
\eqref{Strong_MHD_Maxwell} is obtained by multiplying the equation by
the complex conjugate of a smooth test function $\bb$ with zero
tangential component, integrating the result over $\Dom$, and integrating by parts:
\[
\int_\Dom (\widetilde\mu \bA\SCAL \overline{\bb} + \kappa \ROT
\bA\SCAL\ROT\overline{\bb} )\dif x = \int_\Dom \bbf\SCAL
\overline{\bb}\dif x.
\]
The integral on the left-hand side makes sense if $\bA,\bb\in
\Hrot$. 
% Since the Dirichlet condition $\gamma\upc(\bA)=\bzero$ is enforced on
% $\front\subD$ only, we must properly define  the restriction of the linear
% forms in $\bH^{-\frac12}(\front)$ to functions that are only defined
% on $\front\subD$.  Let $\tbH^{\frac12}(\front\subD)$ be composed of
% the functions $\btheta$ defined on $\front\subD$ whose zero-extension
% to $\front$, say $\tilde\btheta$, is in $\bH^{\frac12}(\front)$. Then,
% for all $\bb \in \Hrot$, the restriction
% $\gamma\upc(\bb)_{|\front\subD}$ is defined in
% $\tbH^{\frac12}(\front\subD)'$ as follows:
% \begin{equation}
% \langle\gamma\upc(\bb)_{|\front\subD},\btheta\rangle_{\front\subD} :=
% \langle\gamma\upc(\bb),\tilde \btheta\rangle_{\front}, \qquad \forall
% \btheta\in \tbH^{\frac12}(\front\subD). 
% \end{equation}
We introduce the following closed subspace of $\Hrot$ to account for the boundary conditions:
\begin{equation}
\bV_0:= \Hrotz = \{\bb\in \Hrot\st\gamma\upc(\bb) =\bzero\}.
\end{equation}
% $\bV_0$ is a closed subspace of $\Hrot$. Indeed, let
% $(\bb_n)_{n\in\polN}$ be a Cauchy sequence in $\bV_0$. Then $\bb_n\to
% \bb$ in $\Hrot$ and, for all $\btheta\in \tbH^{\frac12}(\front\subD)$,
% we have
% \[
% 0 = \langle \gamma\upc(\bb_n)_{|\front\subD},\btheta\rangle_{\front\subD}
% = \langle\gamma\upc(\bb_n),\tilde \btheta\rangle_{\front} \to 
% \langle\gamma\upc(\bb),\tilde \btheta\rangle_{\front},
% \]
% which proves that $\bb\in \bV_0$ since $\gamma\upc:\Hrot\to
% \bH^{-\frac12}(\front)$ is bounded.
The weak formulation of \eqref{Strong_MHD_Maxwell} is the following:
\begin{equation}
\left\{
\begin{array}{l}
\text{Find $\bA\in \bV_0$ such that}\\[2pt]
a_{\widetilde\mu,\kappa}(\bA,\bb) =  \form(\bb),\quad \forall\bb\in \bV_0,
\end{array}
\right.
\label{Weak_MHD_Maxwell}
\end{equation}
where the sesquilinear form $a_{\widetilde\mu,\kappa}$ and the antilinear form $\form$ are defined as follows:
\begin{subequations} \begin{align}a_{\widetilde\mu,\kappa}(\ba,\bb)&:=\int_\Dom (\widetilde\mu \ba\SCAL \overline{\bb} +
\kappa \ROT \ba\SCAL\ROT\overline{\bb} )\dif x, \label{eq:def_akm} \\ 
\form(\bb) &:=\int_\Dom \bbf\SCAL \overline{\bb}\dif x. 
\end{align} \end{subequations}

\begin{Th}[Well-posedness]  \label{Th:wellposedness_Weak_MHD_Maxwell} 
Assume that $\bef\in
\bL^2(\Dom)$, $\widetilde\mu, \kappa \in
L^\infty(\Dom;\polC)$, and that~\eqref{Hyp:mu_kappa_Weak_MHD_Maxwell} holds. 
Then the sesquilinear form $a_{\widetilde\mu,\kappa}$ is coercive:
\begin{equation}\label{eq:Maxwell_coer_basic}
\Re(\text{e}^{i\theta}a_{\widetilde\mu,\kappa}(\bb,\bb)) \ge \min(\mu_\flat,\ell_\Dom^{-2}\kappa_\flat)
  \|\bb\|_{\Hrot}^2,
\end{equation}
for all $\bb\in \Hrot$, and the problem~\eqref{Weak_MHD_Maxwell} is well-posed.
\end{Th}

\begin{proof}
We apply the complex version of the Lax--Milgram Lemma.  The
coercivity follows from~\eqref{Hyp:mu_kappa_Weak_MHD_Maxwell}, and
the boundedness of $a_{\tmu,\kappa}$ and $\form$ are consequences of
$\widetilde\mu, \kappa \in L^\infty(\Dom;\polC)$ and $\bef\in
\bL^2(\Dom)$.  
\eproof

Let us now briefly recall
elementary approximation results solely based on the coercivity of the sesquilinear form $a_{\widetilde\mu,\kappa}$.  
We consider a shape-regular sequence of affine meshes $\famTh$ of
$\Dom$ and the associated approximation setting introduced in
\S\ref{Sec:Generic_FE_setting}-\S\ref{Sec:Best_approximation_result}.
In the rest of the paper, $k$ is the largest integer such that 
$\bpolP_{k}(\Real^3;\Real^3)\subset \wbP\upc\cap \wbP\upd$.
The finite element space we are going to work with is defined by
\begin{equation}
\bV_{h0}=\bP\upc_0(\calT_h) = \{\bb_h\in \bP\upc(\calT_h)\st 
\bb_h\CROSS \bn_{|\front}=\bzero\}.
\label{def_of_Xh_Maxwell}
\end{equation}
Observe that the Dirichlet boundary condition is strongly enforced.
The discrete counterpart of \eqref{Weak_MHD_Maxwell} is formulated as follows:
\begin{equation}
\left\{
\begin{array}{l}
\text{Find $\bA_h\in \bV_{h0}$ such that}\\[2pt]
a_{\widetilde\mu,\kappa}(\bA_h,\bb_h) =  \form(\bb_h),\quad \forall\bb_h\in \bV_{h0}.
\end{array}
\right.
\label{Discrete_Weak_MHD_Maxwell}
\end{equation}
Theorem~\ref{Th:wellposedness_Weak_MHD_Maxwell} together with
the Lax--Milgram Lemma and the conformity property
$\bV_{h0}\subset\bV_0$ implies that \eqref{Discrete_Weak_MHD_Maxwell} has a
unique solution.

\begin{Th}[Error estimate] \label{Th:Convergence_Discrete_Weak_MHD_Maxwell} Under the
  assumptions of Theorem~\textup{\ref{Th:wellposedness_Weak_MHD_Maxwell}}, the
  following holds true:
\begin{align}
\|\bA-\bA_h\|_{\Hrot} &\le \frac{\max(\mu_\sharp,\ell_\Dom^{-2}\kappa_\sharp)}{\min(\mu_\flat,\ell_\Dom^{-2}\kappa_\flat)} \inf_{\bb_h \in \bV_{h0}} \|\bA-\bb_h\|_{\Hrot}.
\label{Hrot_estimate:Th:Convergence_Discrete_Weak_MHD_Maxwell}
\end{align}
Moreover, assuming that there is $r\in (0,k+1]$ so that $\bA\in
  \bH^r(\Dom)$ and $\ROT\bA\in \bH^r(\Dom)$, the following error
estimate holds true:
\begin{align}
  \|\bA-\bA_h\|_{\Hrot} &\le c \, h^{r} (|\bA|_{\bH^r(\Dom)} + \ell_\Dom |\ROT\bA|_{\bH^r(\Dom)}),
\label{Hrot_estimate_hr_Th:Convergence_Discrete_Weak_MHD_Maxwell}
\end{align}
where the constant $c$ is uniform w.r.t.~$h$ and proportional to $\frac{\max(\mu_\sharp,\ell_\Dom^{-2}\kappa_\sharp)}{\min(\mu_\flat,\ell_\Dom^{-2}\kappa_\flat)}$.
\end{Th}
\begin{proof}
The sesquilinear form $a_{\tmu,\kappa}$ is bounded on
$\bV_0\times\bV_0$ with
\begin{equation} \label{eq:a_bounded}
|a_{\tmu,\kappa}(\ba,\bb)| \le
\max(\mu_\sharp,\ell_\Dom^{-2}\kappa_\sharp)\|\ba\|_{\Hrot}
\|\bb\|_{\Hrot},
\end{equation} 
and it is coercive on $\bV_{h0}$ with coercivity constant
$\min(\mu_\flat,\ell_\Dom^{-2}\kappa_\flat)$. Hence, using the
abstract error estimate from \citet[Thm.~2]{XuZik:03}, we obtain
\eqref{Hrot_estimate:Th:Convergence_Discrete_Weak_MHD_Maxwell}.
The second
inequality~\eqref{Hrot_estimate_hr_Th:Convergence_Discrete_Weak_MHD_Maxwell}
is a consequence of Theorem~\ref{Th:Quasi_interpolation_commutes}
together with the estimates of the best approximation error in
Theorem~\ref{Th:glob_best_app} and
Theorem~\ref{Th:best_approx_Dir}. More precisely, we estimate  from above the
infimum
in~\eqref{Hrot_estimate:Th:Convergence_Discrete_Weak_MHD_Maxwell}
 by taking $\bb_h=\calJ_{h0}\upc(\bA)\in \bV_{h0}$. Using
the notation from
\S\ref{Sec:Generic_FE_setting}-\S\ref{Sec:Best_approximation_result},
we have $\bV_{h0} =\bP_0\upc(\calT_h)$ and
\begin{align*}
\|\bA-\calJ_{h0}\upc(\bA)\|_{\Hrot} &\le \|\bA-\calJ_{h0}\upc(\bA)\|_{\Ldeuxd} + \ell_\Dom
\|\ROT(\bA-\calJ_{h0}\upc(\bA))\|_{\Ldeuxd}\\
&=  \|\bA-\calJ_{h0}\upc(\bA)\|_{\Ldeuxd} + \ell_\Dom \|\ROT\bA-\calJ_{h0}\upd(\ROT\bA)\|_{\Ldeuxd} \\
&\le c \inf_{\bb_h \in \bP_{0}\upc(\calT_h)} \|\bA- \bb_h\|_{\Ldeuxd} + 
c'\ell_\Dom \inf_{\bd_h \in \bP_{0}\upd(\calT_h)} \|\ROT\bA - \bd_h\|_{\Ldeuxd} \\
& \le c\, h^{r} (|\bA|_{\bH^r(\Dom)} + \ell_\Dom |\ROT \bA|_{\bH^r(\Dom)}).
\end{align*}
This completes the proof.
%The proof for $|\front\subD|=0$ is similar.  
\end{proof}

\begin{Rem}[Convergence rate]
\cmag{An alternative proof of
\eqref{Hrot_estimate_hr_Th:Convergence_Discrete_Weak_MHD_Maxwell}
is given in \citet[Prop.~4]{Ciarlet:16}
using subtle decompositions of the subspaces
$\bX_{0\tmu}$ and $\bX_{*\kappa^{-1}}$ defined by~\eqref{eq:Hrotz_div0}
and~\eqref{eq:Hrot_divz0} below. The main idea is that
fields in these spaces can be decomposed into a regular part which
can be approximated using the N\'ed\'elec interpolation operator
and a singular part that crucially takes the form of the gradient
of some potential that can be approximated using  
the Scott--Zhang quasi-interpolation operator.}
\end{Rem}

Finally, let us see whether an improved estimate on
$\|\bA-\bA_h\|_{\bL^2(\Dom)}$ can be obtained by the Aubin--Nitsche duality
argument. It is at this point that we realize
that the approach we have taken so far is too simplistic. To better
understand the problem,  let us recall 
a fundamental result that relates the Aubin--Nitsche Lemma
to the compactness of the embedding $V\hookrightarrow L$ where $L$ is the pivot space.
\begin{Th}[Aubin--Nitsche=compactness] 
\label{Thm:Nitsche_Aubin_compactness}
Let $V\hookrightarrow L$ be two Hilbert spaces with continuous
embedding. Let $a:V\CROSS V\to \polC$ be a continuous and coercive sesquilinear
from. Let $(V_h)_{h>0}$ be a sequence of conforming
finite-dimensional approximation spaces.  Let
$G_h: V\to V_h\subset L$  be the discrete solution map defined by $a(G_h(v)-v,v_h)=0$ for all $v_h\in V_h$. Then,
$\lim_{h\to 0} \left(\sup_{v\in V\setminus V_h}\frac{\|G_h(v)-v\|_{L}}{\|G_h(v)-v\|_{V}}\right)=0$
iff the embedding $V\hookrightarrow L$ is compact.
\end{Th}

\begin{proof} See \citet[Thm.~1.1]{Sayas_2004}.
\eproof

\forceindent Let us illustrate this result in the present setting. We
have $\bV_0 =\Hrotz$ and $\bL=\bL^2(\Dom)$, and
Theorem~\ref{Thm:Nitsche_Aubin_compactness} tells us that the
Aubin--Nitsche argument provides an extra rate of convergence in the
$\bL^2$-norm if and only if the embedding
$\Hrotz \hookrightarrow \bL^2(\Dom)$ is compact, which is false. The
conclusion of this argumentation is that we should try to find a space
smaller than $\Hrotz$ where $\bA$ lives and that embeds compactly into
$\bL^2(\Dom)$. We are going to see in Lemma~\ref{lem:pickup}
that a good candidate is essentially $\Hrotz\cap\Hdiv$, \cmag{as pointed
out in \citet[Thm.~2.1-2.3]{Weber:80}}.

\subsection{Regularity pickup}
\label{sec:Maxwell_regularity}

We now recall some key regularity results related to the
curl operator and we use them to infer some regularity pickup in the
scale of Sobolev spaces for $\bA$ and $\ROT\bA$ where 
$\bA$ is the unique solution to \eqref{Strong_MHD_Maxwell}.
 
\forceindent%
Recall that $\bV_0=\Hrotz$ and consider the subspace
$\bX_{0\tmu}:=\{\bb \in \bV_0 \tq \DIV(\tmu\bb)=0\}$. Setting
\begin{equation}
M_0:= \Hunz,
\end{equation} 
a distribution argument shows that 
we can equivalently define $\bX_{0\tmu}$ by
\begin{equation} \label{eq:Hrotz_div0}\textstyle
\bX_{0\tmu}=\{\bb \in \bV_0 \tq (\tmu\bb,\GRAD m)_{\Ldeuxd}=0, \; \forall m\in M_0\},
\end{equation}
where $(\cdot,\cdot)_{\Ldeuxd}$ denotes the inner product in $\Ldeuxd$.
Let us also set
\begin{equation} \label{eq:def_Mstar}
M_*:=\{q\in \Hun\st (q,1)_{\Ldeux}=0\},
\end{equation} 
and define the following subspace:
\begin{equation} \label{eq:Hrot_divz0}
\bX_{*\kappa^{-1}}=\{\bb\in\Hrot\st (\kappa^{-1}\bb,\GRAD m)_{\Ldeuxd}=0, \; \forall m\in M_*\}. 
\end{equation}

\begin{Lem}[Regularity pickup] \label{lem:pickup}
Let $\Dom$ be an open, bounded, and connected Lipschitz subset in
$\Real^3$.
\textup{(i)} Assume that the boundary $\front$ is connected and that 
$\tmu$ is piecewise smooth as
specified in~\textup{\S\ref{Sec:Weak_MHD_Maxwell}}. 
Then there exist 
$s>0$ and $\check C_\Dom>0$ (depending on $\Dom$ and the
heterogeneity ratio $\mu_{\sharp/\flat}$
but not on $\mu_\flat$ alone)
such that
\begin{equation}
\check C_\Dom \ell_\Dom^{-1}\|\bb\|_{\bH^s(\Dom)} \le \|\ROT \bb\|_{\Ldeuxd},
\qquad \forall \bb\in \bX_{0\tmu}.
\label{HDrot_HNdiv_in_Hs}
\end{equation}
\textup{(ii)} Assume that $\Dom$ is simply connected and that $\kappa$ 
is piecewise smooth as specified
in~\textup{\S\ref{Sec:Weak_MHD_Maxwell}}. 
Then there exist $s'>0$ and $\check C_\Dom'>0$
(depending on $\Dom$ and the heterogeneity ratio
$\kappa_{\sharp/\flat}$ but not on $\kappa_\flat$ alone) such that
\begin{equation}
\check C_\Dom' \ell_\Dom^{-1}\|\bb\|_{\bH^{s'}(\Dom)} \le \|\ROT \bb\|_{\Ldeuxd},
\qquad \forall \bb\in\bX_{*\kappa^{-1}}.
\label{HDrot_HNdiv_in_Hs_star}
\end{equation}
\end{Lem}

\forceindent%
The above inequalities, proved in \citet{Jochmann_maxwell_1999} and
\citet{Bonito_guermond_Luddens_amd2_2013}, generalize the classical
results due to \citet[Thm.~3.1]{Birman_Solomyak_1987} and
\citet[Thm.~2]{Costabel_1990}, where the material properties are
assumed to be either constant or smooth and in this case the
smoothness index is $s=\frac12$.  One even has $s\in(\frac12,1]$ if
$\Dom$ is a Lipschitz polyhedron (see \citet[Prop.~3.7]{AmBDG:98}) and
$s=1$ if $\Dom$ is convex (see \citet[Thm.~2.17]{AmBDG:98}). \cmag{We
also refer the reader to
\citet[Thm.~16]{Ciarlet:16} for particular situations in 
heterogeneous materials for which local smoothness with index $s>\frac12$
can be established.}

\forceindent%
Let us now examine the consequences of Lemma~\ref{lem:pickup} on the
Sobolev regularity of $\bA$ and $\ROT\bA$ where 
$\bA$ is the unique solution to \eqref{Strong_MHD_Maxwell}.
Recalling~\eqref{eq:divA}, we infer that $\bA$ is actually a member of
$\bX_{0\tmu}$. Owing to~\eqref{HDrot_HNdiv_in_Hs}, we infer that there
is $s>0$ so that
\begin{equation} \label{eq:A_Hs}
\bA \in \bH^s(\Dom).
\end{equation}
Furthermore, the field $\bR:=\kappa\ROT\bA$ is in $\bX_{*\kappa^{-1}}$,
so that we deduce from~\eqref{HDrot_HNdiv_in_Hs_star} that there is
$s'>0$ so that $\bR \in \bH^{s'}(\Dom)$. In addition, the material
property $\kappa$ being piecewise smooth, we infer that the 
following multiplier property holds (see
\citet[Lem.~2]{Jochmann_maxwell_1999} and
\citet[Prop.~2.1]{Bonito_guermond_Luddens_amd2_2013}): There exists
$\tau>0$ and $C_{\kappa^{-1}}$ such that
\begin{equation} \label{eq:multiplier}
|\kappa^{-1}\bxi|_{\bH^{\tau'}(\Dom)} \le C_{\kappa^{-1}}|\bxi|_{\bH^{\tau'}(\Dom)},
\qquad \forall \bxi\in \bH^{\tau}(\Dom), \quad \forall \tau'\in [0,\tau].
\end{equation}
Letting $s'':=\min(s',\tau)>0$, we conclude that 
\begin{equation} \label{eq:rotA_Hs}
\ROT \bA \in \bH^{s''}(\Dom).
\end{equation}

\section{Coercivity revisited} 
\label{Sec:coercivity:revisited} 
To cope with the loss of coercivity of the sesquilinear form
$a_{\widetilde\mu,\kappa}:\bV_{h0}\CROSS \bV_{h0}\to \polC$ when the
lower bound on $\tmu$ becomes very small, we derive in this section
a sharper coercivity property  in a
proper subspace of $\bV_{h0}$.
The loss of coercivity occurs in the following situations: (i)
In the low frequency limit $(\omega \to 0)$ when
$\widetilde{\mu}=i\omega \mu$, as in the eddy current problem; (ii)
If $\kappa\in \Real$ and $\sigma \ll \omega \epsilon$ when
$\widetilde{\mu} = -\omega^2 \epsilon + i\omega\sigma$, as in the
Helmholtz problem.  

\subsection{Continuous Poincar\'e--Steklov inequality}
\label{Sec:PS_cont} 

Recall the subspace $\bX_{0\tmu}$ of $\bV_0=\Hrotz$ defined in~\eqref{eq:Hrotz_div0}.

\begin{Lem}[Helmholtz decomposition] \label{Lem:Helmholtz_decomposition} 
The following direct sum holds true:
\begin{equation}
\bV_0 = \bX_{0\tmu} 
\oplus \GRAD M_0. \label{Eq:Lem:Helmholtz_decomposition}
\end{equation}
\end{Lem}
\begin{proof}
Let $\bb\in \bV_0$ and let $p\in M_0$ be the unique solution to the
following problem:
$\int_{\Dom} \tmu \GRAD p\SCAL\GRAD \overline{q}\dif x= \int_{\Dom}
\tmu \bb\SCAL \GRAD \overline{q}\dif x$
for all $q\in M_0$. The assumptions on $\tmu$ indeed imply that
there is a unique solution to this problem. Then we set
$\bv=\bb-\GRAD p$ and observe that $\bv\in \bX_{0\tmu}$.  The
sum is direct because if $\bv+\GRAD p=\bzero$, then
$\int_{\Dom} \tmu \GRAD p\SCAL\GRAD \overline{p}\dif x=0$ owing
to~\eqref{eq:Hrotz_div0}, which in tun implies that $p=0$ and
$\bv=\bzero$.
\end{proof}

\begin{Lem}[Poincar\'e--Steklov] \label{lem:PS}
Assume that the boundary $\front$ is connected and that 
$\tmu$ is piecewise smooth as
specified in~\textup{\S\ref{Sec:Weak_MHD_Maxwell}}. 
Then the following Poincar\'e--Steklov inequality holds:
\begin{equation} \label{eq:PoinK_curl} 
\cPoinc{\Dom}
\ell_\Dom^{-1}\|\bb\|_{\Ldeuxd} \le \|\ROT \bb\|_{\Ldeuxd}, \qquad
\forall \bb\in \bX_{0\tmu}, 
\end{equation} 
where $\cPoinc{\Dom}$ can depend on $\Dom$ and the
heterogeneity ratio $\mu_{\sharp/\flat}$
but not on $\mu_\flat$ alone.
\end{Lem}

\bproof
This is a direct consequence of~\eqref{HDrot_HNdiv_in_Hs}.
\eproof

The bound~\eqref{eq:PoinK_curl} is the coercivity property that we
need.  Indeed, \eqref{eq:PoinK_curl} implies
the following series of inequalities for all
$\bb \in \bX_{0\tmu}$:
\begin{align}
\Re(e^{i\theta} a_{\tmu,\kappa}(\bb,\bb)) & \ge \mu_\flat\|\bb\|_{\Ldeuxd}^2 +
\kappa_\flat \|\ROT\bb\|_{\Ldeuxd}^2 \ge \kappa_\flat \|\ROT\bb\|_{\Ldeuxd}^2 \nonumber\\
&\ge 
\frac12\kappa_\flat( \|\ROT\bb\|_{\Ldeuxd}^2 + \cPoinc{\Dom}^2\ell_\Dom^{-2}\|\bb\|_{\Ldeuxd}^2)\nonumber \\
& \ge \frac12\kappa_\flat\ell_\Dom^{-2}\min(1,\cPoinc{\Dom}^2) \|\bb\|_{\Hrot}^2.\label{eq:Maxwell_coer_PoinK}
\end{align}
This shows that the sesquilinear form $a_{\tmu,\kappa}$ is coercive on
$\bX_{0\tmu}$ with a parameter depending on the
heterogeneity ratio $\mu_{\sharp/\flat}$ but not on $\mu_\flat$ alone
(whereas the coercivity parameter on $\bV_0$ is
$\min(\mu_\flat,\ell_\Dom^{-2}\kappa_\flat)$,
see~\eqref{eq:Maxwell_coer_basic}).

\subsection{Discrete Poincar\'e--Steklov inequality}

We show in this section that the ideas of~\S\ref{Sec:PS_cont} can be
reproduced at the discrete level when working with $\Hrt$-conforming
finite elements.  We consider again the discrete
problem~\eqref{Discrete_Weak_MHD_Maxwell}.

Our first step is to realize a discrete counterpart
of~\eqref{Eq:Lem:Helmholtz_decomposition} 
in order to weakly control the divergence of the
discrete vector fields in $\bV_{h0}$. Let us introduce the
$H^1_0$-conforming space
\begin{equation}
  M_{h0} := P_{0}\upg(\calT_h) = \{q_h\in \Hunz \st \mapKg(q_{|K}) \in \wP\upg,\; \forall K\in\calT_h\}. 
\end{equation}
Note that the commutative diagram~\eqref{Diag:wPupg_to_wPupb} implies
that $\GRAD M_{h0}\subset \bV_{h0}$, \ie the polynomial degrees of the
approximations in $M_{h0}$ and $\bV_{h0}$ are compatible. In practice,
the polynomial degree of the reference finite element $(\wK,\wP\upg,\wSigma\upg)$ is $(k+1)$, \ie $\polP_{k+1}(\Real^3;\Real)\subset \wP$.

Since it is not reasonable to consider the space
$\{\bb_h \in \bV_{h0}\st \DIV (\tmu\bb_h)=0\}$ because the normal component
of $\tmu \bb_h$ may jump across the mesh interfaces, we are going to consider
instead the space
\begin{equation}\textstyle
\bX_{h0\tmu} := \{ \bb_h\in \bV_{h0} \tq (\tmu\bb_h, \GRAD m_h)_{\Ldeuxd}= 0,\; \forall m_h\in M_{h0}\}.
\label{discrete_Helmholtz_Maxwell_robust}
\end{equation}
The subtlety here is that $\bX_{h0\tmu}$ is not a subspace of
$\bX_{0\tmu}$, \ie the approximation setting is nonconforming.

\begin{Lem}[Discrete Helmholtz decomposition]
The following direct sum holds true: 
\begin{equation}
\bV_{h0} = \bX_{h0\tmu}\oplus \GRAD M_{h0}.
\end{equation}
\end{Lem}

\bproof
The proof is similar to that of Lemma~\ref{Lem:Helmholtz_decomposition}  since
$\GRAD M_{h0}\subset \bV_{h0}$.
\eproof

\begin{Lem}[Discrete solution] \label{Lem:Ah_is_in_Xh0}
Let $\bA_h\in\bV_{h0}$ be the unique solution to~\eqref{Discrete_Weak_MHD_Maxwell}. Then, $\bA_h\in \bX_{h0\tmu}$. 
\end{Lem}

\bproof We must show that
$\int_\Dom \tmu\bA_h\SCAL \GRAD m_h\dif x=0$ for all $m_h\in M_{h0}$.
Since $\GRAD m_h \in \GRAD M_{h0} \subset \bV_{h0}$,  $\GRAD m_h$ is an admissible
test function in~\eqref{Discrete_Weak_MHD_Maxwell}. Recalling that
$\DIV\bef=0$, we infer that
$0=\ell(\GRAD m_h) = a_{\tmu,\kappa}(\bA_h,\GRAD m_h) = \int_\Dom\tmu
\bA_h\SCAL \GRAD m_h\dif x$
since $\ROT(\GRAD m_h)=\bzero$. This completes the proof.  \eproof

We now establish a discrete counterpart to the Poincar\'e--Steklov
inequality~\eqref{eq:PoinK_curl}. This result is not straightforward
since $\bX_{h0\tmu}$ is not a subspace of $\bX_{0\tmu}$. The key tools that we
are going to invoke are the commuting quasi-interpolation projectors from
Theorem~\ref{Th:Quasi_interpolation_commutes}.

\begin{Th}[Discrete Poincar\'e--Steklov] \label{Thm:Poincare} 
Under the assumptions of Lemma~\textup{\ref{lem:PS}}, there is a
  uniform constant $\cPoinc{\calT_h}>0$ (depending on $\cPoinc{\Dom}$,
  the polynomial degree $k$, the shape-regularity of $\calT_h$, and
  the heterogeneity ratio $\mu_{\sharp/\flat}$, but not on $\mu_\flat$
  alone) such that
\begin{equation}
  \cPoinc{\calT_h} \ell_\Dom^{-1} \|\bx_h\|_{\bL^2(\Dom)} \le \|\ROT \bx_h\|_{\bL^2(\Dom)},
  \qquad \forall \bx_h\in \bX_{h0\tmu}.
\label{Poincare_Maxwell_robust} 
\end{equation}
\end{Th}
\bproof
Let $\bx_h\in\bX_{h0\tmu}$ be a nonzero discrete field.
Let $\phi(\bx_h)\in M_0=\Hunz$ be the solution to the following well-posed Poisson problem:
\[
(\tmu\GRAD \phi(\bx_h),\GRAD m)_{\Ldeuxd} = 
(\tmu\bx_h,\GRAD m)_{\Ldeuxd},\qquad \forall m\in M_0.
\]
Let us define
$\bxi(\bx_h):= \bx_h - \GRAD\phi(\bx_h)$.  This definition implies that
$\bxi(\bx_h)\in \bX_{0\tmu}$. Upon invoking the quasi-interpolation operators
$\calJ_{h0}\upc$ and $\calJ_{h0}\upd$ introduced 
in Theorem~\ref{Th:Quasi_interpolation_commutes} ,
we now observe that
\begin{equation} \label{eq:x-Jxi=grad}
\bx_h - \calJ_{h0}\upc(\bxi(\bx_h)) = \calJ_{h0}\upc(\bx_h - \bxi(\bx_h))
= \calJ_{h0}\upc(\GRAD(\phi(\bx_h))) = \GRAD(\calJ_{h0}\upg(\phi(\bx_h))),
\end{equation}
where we have used that $\calJ_{h0}\upc(\bx_h) = \bx_h$ and 
the commuting properties of the operators $\calJ_{h0}\upg$ and $\calJ_{h0}\upc$. Since $\bx_h\in\bX_{h0\tmu}$,
we infer that
$(\tmu\bx_h,\GRAD(\calJ_{h0}\upg(\phi(\bx_h))))_{\Ldeuxd}=0$, so that
\begin{align*}
(\tmu\bx_h, \bx_h)_{\Ldeuxd} &= 
(\tmu\bx_h,\bx_h - \calJ_{h0}\upc(\bxi(\bx_h)))_{\Ldeuxd}
+( \tmu\bx_h,\calJ_{h0}\upc(\bxi(\bx_h)))_{\Ldeuxd}\\
&= (\tmu\bx_h,\calJ_{h0}\upc(\bxi(\bx_h)))_{\Ldeuxd}.
\end{align*}
Multiplying by $\text{e}^{i\theta}$, taking the real part, and using the 
Cauchy--Schwarz inequality, we infer that
\begin{align*}
\mu_\flat \|\bx_h\|_{\bL^2(\Dom)}^2 &\le  
\mu_\sharp \|\bx_h\|_{\bL^2(\Dom)} \|\calJ_{h0}\upc(\bxi(\bx_h))\|_{\bL^2(\Dom)}.
\end{align*}
The uniform boundedness of $\calJ_{h0}\upc$ on $\bL^2(\Dom)$ together
with the Poincar\'e--Steklov inequality~\eqref{eq:PoinK_curl} implies that
\[
\|\calJ_{h0}\upc(\bxi(\bx_h))\|_{\bL^2(\Dom)} \le \|\calJ_{h0}\upc\|_{\calL(\bL^2;\bL^2)} \|\bxi(\bx_h) \|_{\bL^2(\Dom)} 
\le \|\calJ_{h0}\upc\|_{\calL(\bL^2;\bL^2)}  \cPoinc{\Dom}^{-1} \ell_\Dom \|\ROT\bx_h\|_{\bL^2(\Dom)},
\] 
so that 
\eqref{Poincare_Maxwell_robust} holds with 
$\cPoinc{\calT_h} = \mu_{\sharp/\flat}^{-1} \|\calJ_{h0}\upc\|_{\calL(\bL^2;\bL^2)}^{-1} \cPoinc{\Dom}$. 
\eproof

\bRem[Alternative proofs] There are many ways to prove the discrete
Poincar\'e--Steklov inequality \eqref{Poincare_Maxwell_robust}.  One
route described in \citet[\S4.2]{Hiptmair_acta_numer_2002} hinges on
(subtle) regularity estimates from
\citet[Lemma~4.7]{AmBDG:98}. Another one, which avoids invoking
regularity estimates, is based on the so-called discrete compactness
argument of \citet{Kikuchi_1989} and further developed by
\citet{Monk_Dem_2001} and \citet{Caorsi_Fer_Raf_2000}. The proof based
on the discrete compactness argument is not constructive but relies
instead on an argument by contradiction. The technique used in the
proof of Theorem~\ref{Thm:Poincare}, inspired
\cmag{from~\citet[Thm.~5.11]{ArnFW:06}} and~\citet[Thm~3.6]{ArnFW:10},   
relies on the existence of
the stable commuting quasi-interpolation projectors $\calJ_h\upc$ and
$\calJ_{h0}\upc$. It was observed in \citet{Boffi:01} that the
existence of commuting quasi-interpolation operators within the
discrete de Rham complex would imply the discrete compactness
property.  \eRem

\subsection{Error analysis in the $\bm{H}(\textup{curl})$-norm}
\label{Sec:Hcurl_estimate_revisited}
We are now in a position to revisit the error analysis
of~\S\ref{Sec:basic}. Let us first show that $\bX_{h0\tmu}$ has the same approximation
properties as $\bV_{h0}$ in $\bX_{0\tmu}$.
\begin{Lem}[Approximation in $\bX_{h0\tmu}$]\label{lem:approx_X_V}
The following holds true:
\begin{equation} \label{eq:approx_X_V}
\inf_{\bx_h\in \bX_{h0\tmu}} \|\bA-\bx_h\|_{\Hrot} \le c\, \mu_{\sharp/\flat}
\inf_{\bb_h\in \bV_{h0}} \|\bA-\bb_h\|_{\Hrot}, \quad \forall \bA\in \bX_{0\tmu},
\end{equation}
where the constant $c$ is uniform with respect to $h$ and the model
parameters.
\end{Lem}
\bproof 
Let $\bA\in \bX_{0\tmu}$, and let $p_h\in M_{h0}$ be
the unique solution to the following discrete Poisson problem:
$(\tmu \GRAD p_h,\GRAD q_h)_{\Ldeuxd} = (\tmu \calJ\upc_{h0}(\bA),
\GRAD q_h)_{\Ldeuxd}$
for all $q_h\in M_{h0}$.  Let us define
$\by_h = \calJ\upc_{h0}(\bA) -\GRAD p_h$. By construction,
$\by_h\in \bX_{h0\tmu}$ and $\ROT\by_h=\ROT\calJ\upc_{h0}(\bA)$. Hence,
$\|\ROT(\bA -\by_h)\|_{\Ldeuxd} = \|\ROT(\bA -
\calJ\upc_{h0}(\bA))\|_{\Ldeuxd}$.
Moreover, since $\DIV(\tmu\bA)=0$, we infer that
\[
(\tmu\GRAD p_h,\GRAD p_h)_{\Ldeuxd} = 
(\tmu \calJ\upc_{h0}(\bA),\GRAD p_h)_{\Ldeuxd} = 
(\tmu (\calJ\upc_{h0}(\bA)-\bA),\GRAD p_h)_{\Ldeuxd},
\]
which in turn implies that
$\|\GRAD p_h\|_{\Ldeuxd} \le \mu_{\sharp/\flat} \|\calJ\upc_{h0}(\bA) -\bA\|_{\Ldeuxd}$. The
above argument shows that
\begin{align*}
\|\bA -\by_h\|_{\Ldeuxd} &\le \|\bA - \calJ\upc_{h0}(\bA)\|_{\Ldeuxd}+
\|\calJ\upc_{h0}(\bA)-\by_h\|_{\Ldeuxd} \\
&\le \|\bA - \calJ\upc_{h0}(\bA)\|_{\Ldeuxd}+
\|\GRAD p_h \|_{\Ldeuxd} \\
&\le (1+\mu_{\sharp/\flat}) \|\bA - \calJ\upc_{h0}(\bA)\|_{\Ldeuxd}.
\end{align*}
In conclusion, we have proved that
\begin{align*}
\inf_{\bx_h\in \bX_{h0\tmu}} \|\bA-\bx_h\|_{\Hrot} &\le \|\bA-\by_h\|_{\Hrot} \\
&\le 
(1+\mu_{\sharp/\flat}) \|\bA- \calJ\upc_{h0}(\bA)\|_{\Hrot}.
\end{align*}
Upon invoking~\eqref{Lp_error_estimate_for_calJh} and the commutative diagrams
\eqref{Diag:calJh0_commutes}, we infer that
\begin{align*}
\|\bA-\calJ\upc_{h0}(\bA)\|_{\Hrot} &\le \|\bA-\calJ_{h0}\upc(\bA)\|_{\Ldeuxd} + \ell_\Dom
\|\ROT(\bA-\calJ_{h0}\upc(\bA))\|_{\Ldeuxd}\\
&=  \|\bA-\calJ_{h0}\upc(\bA)\|_{\Ldeuxd} + \ell_\Dom \|\ROT\bA-\calJ_{h0}\upd(\ROT\bA)\|_{\Ldeuxd} \\
&\le c \inf_{\bb_h \in \bP_{0}\upc(\calT_h)} \|\bA- \bb_h\|_{\Ldeuxd} + 
c'\ell_\Dom \inf_{\bd_h \in \bP_{0}\upd(\calT_h)} \|\ROT\bA - \bd_h\|_{\Ldeuxd} \\
&\le c \inf_{\bb_h \in \bP_{0}\upc(\calT_h)} \|\bA- \bb_h\|_{\Ldeuxd} + 
c'\ell_\Dom \inf_{\bb_h \in \bP_{0}\upc(\calT_h)} \|\ROT(\bA - \bb_h)\|_{\Ldeuxd},
\end{align*}
where the last bound follows by restricting the minimization set to
$\bP_{0}\upc(\calT_h)$ since
$\ROT \bP_{0}\upc(\calT_h)\subset \bP_{0}\upd(\calT_h)$. The
conclusion follows readily.  \eproof
\begin{Th}[Error estimate] \label{Th:Hrot_convergence_MHD_Maxwell_robust}
Assume that the boundary $\front$ is connected and that 
$\tmu$ is piecewise smooth as
specified in~\textup{\S\ref{Sec:Weak_MHD_Maxwell}}. 
Then the following error estimate holds true:
\begin{equation}
\|\bA-\bA_h\|_{\Hrot} \le c\, \max(1,\gamma_{\tmu,\kappa}) 
\inf_{\bb_h\in \bV_{h0}} \|\bA-\bb_h\|_{\Hrot}, \label{Eq:Th:Hrot_convergence_MHD_Maxwell_robus}
\end{equation}
where the constant $c$ is uniform with respect to $h$ and can depend on
the discrete Poincar\'e--Steklov constant $\cPoinc{\calT_h}$ and
the heterogeneity ratios $\mu_{\sharp/\flat}$ and
$\kappa_{\sharp/\flat}$,
and where $\gamma_{\tmu,\kappa}=\mu_\sharp\ell_\Dom^2\kappa_\sharp^{-1}$
is the magnetic Reynolds number.
\end{Th}
\bproof Owing to Lemma~\ref{Lem:Ah_is_in_Xh0}, 
$\bA_h$ solves the following problem: Find
$\bA_h\in \bX_{h0\tmu}$ such that
$a_{\tmu,\kappa}(\bA_h,\bx_h) = \form(\bx_h)$, for all
$\bx_h\in \bX_{h0\tmu}$. Using the discrete Poincar\'e--Steklov
inequality~\eqref{Poincare_Maxwell_robust} and proceeding as
in~\eqref{eq:Maxwell_coer_PoinK}, we infer that
\[
\Re(e^{i\theta} a_{\tmu,\kappa}(\bx_h,\bx_h)) 
\ge \frac12\kappa_\flat\ell_\Dom^{-2}\min(1,\cPoinc{\calT_h}^2) 
\|\bx_h\|_{\Hrot}^2,
\]
for all $\bx_h\in \bX_{h0\tmu}$.  Hence, the above problem is well-posed.
Recalling the boundedness property~\eqref{eq:a_bounded} of the
sesquilinear form $a_{\tmu,\kappa}$ and invoking again the abstract
error estimate from \citet[Thm.~2]{XuZik:03} leads to
\[
\|\bA-\bA_h\|_{\Hrot} \le \frac{2\max(\mu_\sharp,\ell_\Dom^{-2}\kappa_\sharp)}{\kappa_\flat\ell_\Dom^{-2}\min(1,\cPoinc{\calT_h}^2)}\inf_{\bx_h\in \bX_{h0\tmu}} \|\bA-\bx_h\|_{\Hrot}.
\]
We conclude the proof by invoking Lemma~\ref{lem:approx_X_V}.
\eproof

\bRem[Neuman boundary condition] 
The above analysis can be adapted to account for the
Neumann boundary condition
$(\kappa\ROT\bA)_{|\front}\CROSS\bn = \bzero$. This condition implies
that $(\ROT(\kappa\ROT\bA))_{|\front}\SCAL \bn=0$. Moreover,
assuming that $\bef\SCAL\bn_{|\front} =0$, and taking the normal
component of the equation $\tmu \bA +\ROT(\kappa\ROT\bA) = \bef$ at
the boundary gives $\bA\SCAL \bn_{|\front} =0$. Since
$\DIV\bef=0$, we also have that $\DIV(\tmu\bA)=0$.  Using a distribution
argument shows that $\bA\in \bX_{*\tmu}$ where
$\bX_{*\tmu}:=\{\bb\in\Hrot\st (\tmu\bb,\GRAD m)_{\Ldeuxd}=0, \;
\forall m\in M_*\}$ and $M_*$ is defined in~\eqref{eq:def_Mstar}. The discrete
spaces that must be used are now $\bV_{h} = \bP\upc(\calT_h)$ and
$M_{h*}=P\upg(\calT_h)\cap M_*$.  
Using $\bV_{h}$ for the discrete trial and test spaces in the weak
formulation, one then deduces that 
\begin{equation}
\bA_h \in \bX_{h*\tmu} := \{ \bb_h\in \bV_{h} \tq (\tmu\bb_h,
\GRAD m_h)_{\Ldeuxd}=0, \; \forall m_h\in M_{h*}\}.
\end{equation}
The Poincar\'e--Steklov inequality \eqref{Poincare_Maxwell_robust}
still holds if the assumption that $\front$ is connected is replaced
by the assumption that $\Dom$ is simply connected. The error analysis
from Theorem~\ref{Th:Hrot_convergence_MHD_Maxwell_robust} can be
readily adapted.  
\eRem

\bRem[Helmholtz problem] \cmag{The assumption
  \eqref{Hyp:mu_kappa_Weak_MHD_Maxwell} can be replaced by assuming that
  $\essinf_{\bx\in \Dom}|\kappa| \ge \kappa_\flat>0$ and $\tmu$ is not
  an eigenvalue of the operator $\ROT(\kappa\ROT)$
  equipped with the appropriate boundary condition. In this case,
  \eqref{Weak_MHD_Maxwell} is a Helmholtz-type boundary-value
  problem. This problem can be analyzed by using commuting
  quasi-interpolation operators as done in \citet[\S9.1]{ArnFW:06} for
  Neumann boundary conditions, or by using a constructive proof
  relying on Hilbertian bases as done in \citet{Ciarlet:12}. Note that
  the convergence rates derived in \cite{Ciarlet:12} require a
  smoothness index $s>\frac12$ owing to the use of the N\'ed\'elec
  interpolation operator; therefore, the present quasi-interpolation
  operator can be combined with these results to treat more general
  situations regarding the heterogeneity of the material.}  \eRem

\section{The duality argument for edge elements}
\label{Sec:L2}

Our goal in this section is to estimate $(\bA-\bA_h)$ in the
$\bL^2$-norm using a duality argument that invokes a weak control on
the divergence. The subtlety is that, as already mentioned, the
setting is nonconforming since $\bX_{h0\tmu}$ is not a subspace of
$\bX_{0\tmu}$. Recalling Theorem~\ref{Thm:Nitsche_Aubin_compactness},
the compactness that is required from the functional setting to
obtain a better convergence rate in $\bL^2(\Dom)$ will result from
Lemma~\ref{lem:pickup} and the compact embedding
$\bH^s(\Dom)\hookrightarrow\bL^2(\Dom)$, $s>0$.  In this section, we
are going to use both inequalities~\eqref{HDrot_HNdiv_in_Hs}
and~\eqref{HDrot_HNdiv_in_Hs_star} from Lemma~\ref{lem:pickup};
therefore, we assume that the boundary $\front$ is connected and that
$\Dom$ is simply connected, and that both $\tmu$ and $\kappa$ are
piecewise smooth. Recalling the results
of~\S\ref{sec:Maxwell_regularity} with smoothness indices $s,s'>0$ and
the index $\tau>0$ from the multiplier property~\eqref{eq:multiplier}
and letting $s''=\min(s',\tau)$, we have $\bA\in \bH^s(\Dom)$ and
$\ROT\bA\in\bH^{s''}(\Dom)$ with $s,s''>0$. In what follows, we set
\begin{equation}
\sigma := \min(s,s'').
\end{equation}
Recall the magnetic Reynolds number
$\gamma_{\tmu,\kappa}=\mu_\sharp\ell_\Dom^2\kappa_\sharp^{-1}$ and let
us set $\hat\gamma_{\tmu,\kappa}=\max(1,\gamma_{\tmu,\kappa})$.

Let us first start with an approximation result on the curl-preserving
lifting operator $\bxi:\bX_{h0\tmu}\to \bX_{0\tmu}$ defined in the proof of
Theorem~\ref{Thm:Poincare}.  Recall that, for all $\bx_h\in \bX_{h0\tmu}$,
the field $\bxi(\bx_h)\in \bX_{0\tmu}$ is such
$\bxi(\bx_h) = \bx_h - \GRAD\phi(\bx_h)$ with
$\phi(\bx_h)\in \Hunz$; hence $\ROT\bxi(\bx_h) = \ROT\bx_h$.

\begin{Lem}[Curl-preserving lifting] \label{Lem:Rotv_Rotvh_approximation}
Let $s>0$ be the smoothness index introduced in
\eqref{HDrot_HNdiv_in_Hs}. Then, the following holds true:
\begin{equation}
\|\bxi(\bx_h)-\bx_h\|_{\Ldeuxd} \le c\, h^s\ell_\Dom^{1-s}\|\ROT
\bx_h\|_{\Ldeuxd}, 
\qquad \forall \bx_h\in \bX_{h0\tmu},
\label{Eq2:Lem:Rotv_Rotvh_approximation}
\end{equation}
where the constant $c$ is uniform with respect to $h$ and can depend on the
constant $\check C_\Dom$ from~\textup{\eqref{HDrot_HNdiv_in_Hs}} and
the heterogeneity ratio $\mu_{\sharp/\flat}$.
\end{Lem}
\begin{proof} 
Let $\bx_h\in \bX_{h0\tmu}$, and
let us set $\be_h:=\bxi(\bx_h)-\bx_h$. We have seen in the proof of
Theorem~\ref{Thm:Poincare} that
$\calJ_{h0}\upc(\bxi(\bx_h))-\bx_h \in \GRAD M_{h0}$,
see~\eqref{eq:x-Jxi=grad}. This, in turn, implies that
$(\tmu\be_h,\calJ_{h0}\upc(\bxi(\bx_h))-\bx_h)_{\Ldeuxd}=0$ since
$\bxi(\bx_h)\in \bX_{0\tmu}$, $M_{h0}\subset M_0$, and
$\bx_h\in \bX_{h0\tmu}$. Since
$\be_h = (I-\calJ_{h0}\upc)(\bxi(\bx_h)) +
(\calJ_{h0}\upc(\bxi(\bx_h))-\bx_h)$, we infer that
\[
(\tmu \be_h,\be_h)_{\Ldeuxd} = (\tmu \be_h,(I-\calJ_{h0}\upc)(\bxi(\bx_h)))_{\Ldeuxd},
\]
thereby implying that
$\|\be_h\|_{\Ldeuxd} \le \mu_{\sharp/\flat}
\|(I-\calJ_{h0}\upc)(\bxi(\bx_h))\|_{\Ldeuxd}$.
Using the approximation properties of $\calJ_{h0}\upc$ yields
\[
\|\be_h\|_{\Ldeuxd} \le c\, \mu_{\sharp/\flat} h^s |\bxi(\bx_h)|_{\bH^s(\Dom)},
\] 
and we conclude using the bound
$|\bxi(\bx_h)|_{\bH^s(\Dom)}\le \check C_\Dom\ell_\Dom^{1-s} \|\ROT
\bx_h\|_{\Ldeuxd}$
which follows from~\eqref{HDrot_HNdiv_in_Hs}.  
\eproof

\begin{Lem}[Adjoint solution] \label{lem:Maxwell_adjoin_bnd}
Let $\by \in \bX_{0\tmu}$ and let $\bzeta\in \bX_{0\tmu}$ be the (unique) solution to the (adjoint) problem 
$\widetilde{\mu} \bzeta + \ROT (\kappa \ROT \bzeta) = \mu_\flat^{-1}\tmu \by$. Then,
\begin{subequations} \label{eq:bnd_bzeta} \begin{align}
|\bzeta|_{\bH^\sigma(\Dom)} &\le c\, \mu_\sharp^{-1}
  \gamma_{\tmu,\kappa} \ell_\Dom^{-\sigma}\|\by\|_{\Ldeux}, \label{eq:bnd_bzeta1}\\
|\ROT\bzeta|_{\bH^\sigma(\Dom)} &\le  c\, \mu_\sharp^{-1}\gamma_{\tmu,\kappa}\hat\gamma_{\tmu,\kappa}
  \ell_\Dom^{-1-\sigma}\ \|\by\|_{\Ldeux},
\label{eq:bnd_bzeta2}
\end{align}\end{subequations}
where the constant $c$ is uniform with respect to $h$ and can depend on the
constants $\cPoinc{\Dom}$ from~\textup{\eqref{eq:PoinK_curl}},
$\check C_\Dom$, $\check C_\Dom'$
from~\textup{\eqref{HDrot_HNdiv_in_Hs}-\eqref{HDrot_HNdiv_in_Hs_star}},
and the heterogeneity ratios $\mu_{\sharp/\flat}$,
$\kappa_{\sharp/\flat}$, $\kappa_\sharp C_{\kappa^{-1}}$.
\end{Lem}

\bproof Testing the adjoint problem against $\bzeta$, we observe that
$\kappa_{\flat} \|\ROT\bzeta\|_{\Ldeuxd}^2 \le
\mu_{\sharp/\flat} \|\by\|_{\Ldeux} \|\bzeta\|_{\Ldeuxd}$,
so that using the Poincar\'e--Steklov inequality~\eqref{eq:PoinK_curl} to bound
  $\|\bzeta\|_{\Ldeuxd}$ by $\|\ROT\bzeta\|_{\Ldeuxd}$, we infer that
\begin{equation} \label{eq:bnd_rot_zeta}
 \|\ROT\bzeta\|_{\Ldeuxd}
 \le \kappa_{\flat}^{-1} \mu_{\sharp/\flat} \cPoinc{\Dom}^{-1} \ell_\Dom \|\by\|_{\Ldeux}.
\end{equation}
Invoking~\eqref{HDrot_HNdiv_in_Hs} with $\sigma\le s$ yields
\[
|\bzeta|_{\bH^\sigma(\Dom)} \le \check C_\Dom^{-1}\ell_\Dom^{1-\sigma}\|\ROT\bzeta\|_{\Ldeuxd}
\le \kappa_{\flat}^{-1} \mu_{\sharp/\flat} \check C_\Dom^{-1}\cPoinc{\Dom}^{-1} \ell_\Dom^{2-\sigma}\|\by\|_{\Ldeux},
\]
which proves~\eqref{eq:bnd_bzeta1} since $\kappa_{\flat}^{-1}\ell_\Dom^2
= \kappa_{\sharp/\flat} \mu_\sharp^{-1}\gamma_{\tmu,\kappa}$. Let us now
prove~\eqref{eq:bnd_bzeta2}. Invoking~\eqref{HDrot_HNdiv_in_Hs_star}
with $\sigma\le s'$
for $\bb = \kappa\ROT\bzeta$, which is a member of
$\bX_{*\kappa^{-1}}$, we infer that
\[
\check C_\Dom'\ell_\Dom^{-1+\sigma} |\bb|_{\bH^{\sigma}(\Dom)} \le \|\ROT
\bb\|_{\Ldeuxd} = \|\ROT(\kappa\ROT\bzeta)\|_{\Ldeuxd} \le
\mu_{\sharp/\flat}\|\by\|_{\Ldeuxd} + \mu_\sharp
\|\bzeta\|_{\Ldeuxd},
\]
by definition of the adjoint solution
$\bzeta$ and the triangle inequality. 
Invoking again the Poincar\'e--Steklov inequality~\eqref{eq:PoinK_curl} to bound
$\|\bzeta\|_{\Ldeuxd}$ by $\|\ROT\bzeta\|_{\Ldeuxd}$ and
using~\eqref{eq:bnd_rot_zeta} yields
\[
\|\bzeta\|_{\Ldeuxd} \le \kappa_{\flat}^{-1} \mu_{\sharp/\flat} \cPoinc{\Dom}^{-2} \ell_\Dom^2 \|\by\|_{\Ldeux}.
\]
As a result, we obtain
\[
\check C_\Dom'\ell_\Dom^{-1+\sigma} |\bb|_{\bH^{\sigma}(\Dom)} \le
\mu_{\sharp/\flat}(1 + \mu_\sharp\kappa_\flat^{-1} \cPoinc{\Dom}^{-2} \ell_\Dom^2) \|\by\|_{\Ldeuxd},
\]
and we can conclude the proof of~\eqref{eq:bnd_bzeta2} since
$|\ROT \bzeta|_{\bH^\sigma(\Dom)} \le C_{\kappa^{-1}}
|\bb|_{\bH^\sigma(\Dom)}$
owing to the multiplier property~\eqref{eq:multiplier} and $\sigma\le \tau$.  
\eproof
We can now state the main result of this section. 

\begin{Th}[Improved $L^2$-error estimate] \label{Th:improved_L2_estimate}
The following holds true:
\begin{equation}
\|\bA-\bA_h\|_{\bL^2} \le  c \,  \inf_{\bv_h \in \bV_{h0}} \!(\|\bA-\bv_h\|_{\bL^2}
+  \hat\gamma_{\tmu,\kappa}^3 h^\sigma \ell_\Dom^{-\sigma} \|\bA-\bv_h\|_{\Hrt}),
\end{equation}
where the constant $c$ is uniform with respect to $h$ and can depend on the constants $\cPoinc{\Dom}$
  from~\textup{\eqref{eq:PoinK_curl}}, $\check C_\Dom$, $\check
C_\Dom'$ from~\textup{\eqref{HDrot_HNdiv_in_Hs}-\eqref{HDrot_HNdiv_in_Hs_star}}, 
and the heterogeneity ratios
$\mu_{\sharp/\flat}$, $\kappa_{\sharp/\flat}$,
$\kappa_\sharp C_{\kappa^{-1}}$.
\end{Th}
\bproof In this proof, we use the symbol $c$ to denote a generic
positive constant that can have the same parametric dependencies as in
the above statement. Let $\bv_h\in \bX_{h0\tmu}$ and let us set $\bx_h:=\bA_h-\bv_h$;
observe that $\bx_h\in \bX_{h0\tmu}$.  Let $\bzeta\in \bX_{0\tmu}$ be the
solution to the adjoint problem
$\widetilde{\mu} \bzeta + \ROT (\kappa \ROT \bzeta) =
\mu_\flat^{-1}\tmu \bxi(\bx_h)$,
where $\bxi:\bX_{h0\tmu}\to \bX_{0\tmu}$ is the curl-preserving lifting operator considered above.
\\
(1) Let us first estimate $\|\bxi(\bx_h)\|_{\Ldeuxd}$ from above.  Recalling
that $\bxi(\bx_h)-\bx_h=-\GRAD\phi(\bx_h)$ and that
$(\tmu\bxi(\bx_h),\bxi(\bx_h)-\bx_h)_{\Ldeux} =
-(\tmu\bxi(\bx_h),\GRAD\phi(\bx_h))_{\Ldeux}= 0$, we infer that
\begin{align*}
(\bxi(\bx_h),\tmu\bxi(\bx_h))_{\Ldeuxd}  &= (\bx_h,\tmu\bxi(\bx_h))_{\Ldeuxd} \\
& = (\bA-\bv_h,\tmu\bxi(\bx_h))_{\Ldeuxd} + (\bA_h-\bA,\tmu\bxi(\bx_h))_{\Ldeuxd} \\
& = (\bA-\bv_h,\tmu\bxi(\bx_h))_{\Ldeuxd} + \mu_\flat a_{\tmu,\kappa}(\bA_h-\bA,\bzeta) \\
& = (\bA-\bv_h,\tmu\bxi(\bx_h))_{\Ldeuxd} + \mu_\flat a_{\tmu,\kappa}(\bA_h-\bA,\bzeta-\calJ_{h0}\upc(\bzeta)),
\end{align*}
where we have used Galerkin's orthogonality to pass from the third to
the fourth line. Since we have
$|a_{\tmu,\kappa}(\ba,\bb)| \le 
\kappa_\sharp\ell_\Dom^{-2}\hat\gamma_{\tmu,\kappa} \|\ba\|_{\Hrot} \|\bb\|_{\Hrot}$
owing to~\eqref{eq:a_bounded}, we infer from the commutation and
approximation results from
Theorem~\ref{Th:Quasi_interpolation_commutes} and
  Theorem~\ref{Th:best_approx_Dir} that
\begin{align*}
\|\bxi(\bx_h)\|_{\Ldeuxd}^2 \le{}&\mu_{\sharp/\flat} \|\bA-\bv_h\|_{\Ldeuxd}
\|\bxi(\bx_h)\|_{\Ldeuxd} \\&+ c\, 
\kappa_\sharp\ell_\Dom^{-2}\hat\gamma_{\tmu,\kappa} h^\sigma \|\bA-\bA_h\|_{\Hrot}
( |\bzeta|_{\bH^\sigma(\Dom)}
+ \ell_\Dom|\ROT\bzeta|_{\bH^\sigma(\Dom)}).
\end{align*}
Invoking the bounds from Lemma~\ref{lem:Maxwell_adjoin_bnd} on the
adjoint solution with $\by=\bxi(\bx_h)$, we conclude that
\begin{equation}
\|\bxi(\bx_h)\|_{\Ldeuxd} \le \mu_{\sharp/\flat} \|\bA-\bv_h\|_{\Ldeuxd}
+ c\, \hat\gamma_{\tmu,\kappa}^2 h^\sigma\ell_\Dom^{-\sigma}\|\bA-\bA_h\|_{\Hrot}. \label{to_please_ref1}
% \hat\chi_\Dom^{-3}\gamma_{\mu,\kappa}^2\kappa_{\sharp/\flat}'
\end{equation}
(2) The triangle inequality, together with the identity $\bA-\bA_h=\bA-\bv_h -\bx_h$, implies that 
\[
\|\bA-\bA_h\|_{\bL^2} \le \|\bA-\bv_h\|_{\Ldeuxd} + \|\bxi(\bx_h)-\bx_h\|_{\Ldeuxd}
+ \|\bxi(\bx_h)\|_{\Ldeuxd}.
\]
We use Lemma~\ref{Lem:Rotv_Rotvh_approximation} to bound the second
term on the right-hand side as follows:
\begin{align*}
\|\bxi(\bx_h)-\bx_h\|_{\Ldeuxd} 
&\le c\, h^\sigma \ell_\Dom^{1-\sigma}\|\ROT \bx_h\|_{\Ldeuxd} \\ 
&\le c\, h^\sigma \ell_\Dom^{1-\sigma}(\|\ROT(\bA-\bv_h)\|_{\Ldeuxd}+\|\ROT(\bA-\bA_h)\|_{\Ldeuxd}),
\end{align*}
and we use
\eqref{Eq:Th:Hrot_convergence_MHD_Maxwell_robus} to infer that
$\|\bA-\bA_h\|_{\Hrot} \le c \hat\gamma_{\tmu,\kappa} \|\bA-\bv_h\|_{\Hrot}$.
%\mu_{\sharp/\flat} \hat\chi_\Dom^{-2}\gamma_{\mu,\kappa}\kappa_{\sharp/\flat}'
For the third term on the right-hand side, we use the bound on
$\|\bxi(\bx_h)\|_{\Ldeuxd}$ estimated above in Step (1). We conclude by taking the infimum
over $\bv_h\in \bX_{h0\tmu}$, and we use Lemma~\ref{lem:approx_X_V} to extend the
infimum over $\bV_{h0}$.  
\eproof

\bRem[Curl-preserving lifting] The idea of the construction of the
curl-preserving lifting operator invoked in the proof of
Theorem~\ref{Thm:Poincare} and Theorem~\ref{Th:improved_L2_estimate}
is rooted in \citet[p.~249-250]{Monk_1992}.  The statement in
Lemma~\ref{Lem:Rotv_Rotvh_approximation} is similar to that in
\citet[Lem~7.6]{Monk_2003}, but the proof we give is greatly
simplified by our using the commuting quasi-interpolation operators
from Theorem~\ref{Th:Quasi_interpolation_commutes}.  \cmag{The
curl-preserving lifting of $\bA-\bA_h$ is invoked in
\citet[Eq.~(9.9)]{ArnFW:06} and denoted therein by
$\bpsi$. The estimate of $\|\bpsi\|_{\bL^2(\Dom)}$ given one line
above \cite[Eq.~(9.11)]{ArnFW:06} is similar to
\eqref{to_please_ref1} and is obtained by invoking the commuting
quasi-interpolation operators constructed in
\cite[\S5.4]{ArnFW:06} for natural boundary conditions. 
Note that contrary to what is
done in the above reference, we invoke the curl-preserving lifting of
$\bA_h-\bv_h$ instead of $\bA-\bA_h$ and make use of
Lemma~\ref{Lem:Rotv_Rotvh_approximation}, which simplifies the
argument.} Furthermore, the statement 
of Theorem~\ref{Th:improved_L2_estimate} is
similar to that of \citet[Thm.~4.1]{Zhong_Shi_Wittum_Xu_2009}, but
again the proof that we propose is simplified by our using the
commuting quasi-interpolation operators; moreover, the setting
presented in this paper accounts for heterogeneous coefficients
since it holds true for any smoothness index $\sigma$ smaller than
$\frac12$.
% Note finally that the results from
% Theorem~\ref{Th:Quasi_interpolation_commutes},
% Theorem~\ref{Th:glob_best_app}, and Theorem~\ref{Th:best_approx_Dir}
% are essential to overcome the difficulties induced by the lack of
% stability of the canonical interpolation operators in
% $P_0\upg(\calT_h)$, $P\upc_0(\calT_h)$, \etc and to bypass the
% technique known in the literature as the discrete compactness
% argument.  
\eRem

% \ifCMAM
% \begin{acknowledgement}
% This material is based upon work supported in part by the National
% Science Foundation grants DMS 1620058, DMS 1619892, by the Air
% Force Office of Scientific Research, USAF, under grant/contract
% number FA99550-12-0358 and the Army Research Office, under grant number W911NF-15-1-0517. 
% This work was done when the authors were visiting
% the  Institut Henri Poincar\'e during the Fall 2016 Thematic
% Trimester ``Numerical Methods for Partial Differential Equations''. 
% The support of IHP is gratefully acknowledged. 
% \end{acknowledgement}
% \else
\paragraph*{Acknowledgement}
This material is based upon work supported in part by the National
Science Foundation grants DMS 1620058, DMS 1619892, by the Air
Force Office of Scientific Research, USAF, under grant/contract
number FA99550-12-0358 and the Army Research Office, under grant
number W911NF-15-1-0517. This work was done when the authors were visiting
the  Institut Henri Poincar\'e during the Fall 2016 Thematic
Trimester ``Numerical Methods for Partial Differential Equations''. 
The support of IHP is gratefully acknowledged. 
% \fi

\bibliographystyle{abbrvnat}
\bibliography{ref}

\begin{thebibliography}{34}
\providecommand{\natexlab}[1]{#1}
\providecommand{\url}[1]{\texttt{#1}}
\expandafter\ifx\csname urlstyle\endcsname\relax
  \providecommand{\doi}[1]{doi: #1}\else
  \providecommand{\doi}{doi: \begingroup \urlstyle{rm}\Url}\fi

\bibitem[Amrouche et~al.(1998)Amrouche, Bernardi, Dauge, and Girault]{AmBDG:98}
C.~Amrouche, C.~Bernardi, M.~Dauge, and V.~Girault.
\newblock Vector potentials in three-dimensional non-smooth domains.
\newblock \emph{Math. Methods Appl. Sci.}, 21\penalty0 (9):\penalty0 823--864,
  1998.

\bibitem[Arnold et~al.(2006)Arnold, Falk, and Winther]{ArnFW:06}
D.~N. Arnold, R.~S. Falk, and R.~Winther.
\newblock Finite element exterior calculus, homological techniques, and
  applications.
\newblock \emph{Acta Numer.}, 15:\penalty0 1--155, 2006.

\bibitem[Arnold et~al.(2010)Arnold, Falk, and Winther]{ArnFW:10}
D.~N. Arnold, R.~S. Falk, and R.~Winther.
\newblock Finite element exterior calculus: from {H}odge theory to numerical
  stability.
\newblock \emph{Bull. Amer. Math. Soc. (N.S.)}, 47\penalty0 (2):\penalty0
  281--354, 2010.

\bibitem[Birman and Solomyak(1987)]{Birman_Solomyak_1987}
M.~S. Birman and M.~Z. Solomyak.
\newblock $l_2$ -{T}heory of the {M}axwell operator in arbitrary domains.
\newblock \emph{Russian Mathematical Surveys}, 42\penalty0 (6):\penalty0 75,
  1987.

\bibitem[Boffi(2001)]{Boffi:01}
D.~Boffi.
\newblock A note on the de {R}ham complex and a discrete compactness property.
\newblock \emph{Appl. Math. Lett.}, 14\penalty0 (1):\penalty0 33--38, 2001.

\bibitem[Boffi and Gastaldi(2006)]{BofGa:06}
D.~Boffi and L.~Gastaldi.
\newblock Interpolation estimates for edge finite elements and application to
  band gap computation.
\newblock \emph{Appl. Numer. Math.}, 56\penalty0 (10-11):\penalty0 1283--1292,
  2006.

\bibitem[Bonito et~al.(2013)Bonito, Guermond, and
  Luddens]{Bonito_guermond_Luddens_amd2_2013}
A.~Bonito, J.-L. Guermond, and F.~Luddens.
\newblock Regularity of the maxwell equations in heterogeneous media and
  lipschitz domains.
\newblock \emph{J. Math. Anal. Appl.}, 408:\penalty0 498--512, 2013.

\bibitem[Bossavit(1998)]{Bossavit_GB}
A.~Bossavit.
\newblock \emph{Computational Electromagnetism, Variational Formulations,
  Complementary, Edge Elements}, volume~2 of \emph{Electromagnetism}.
\newblock Academic Press, New York, NY, 1998.

\bibitem[Brezis(2011)]{Brezis:11}
H.~Brezis.
\newblock \emph{Functional analysis, {S}obolev spaces and partial differential
  equations}.
\newblock Universitext. Springer, New York, 2011.

\bibitem[Caorsi et~al.(2000)Caorsi, Fernandes, and
  Raffetto]{Caorsi_Fer_Raf_2000}
S.~Caorsi, P.~Fernandes, and M.~Raffetto.
\newblock On the convergence of {G}alerkin finite element approximations of
  electromagnetic eigenproblems.
\newblock \emph{SIAM J. Numer. Anal.}, 38\penalty0 (2):\penalty0 580--607
  (electronic), 2000.

\bibitem[Christiansen(2007)]{Christ:07}
S.~H. Christiansen.
\newblock Stability of {H}odge decompositions in finite element spaces of
  differential forms in arbitrary dimension.
\newblock \emph{Numer. Math.}, 107\penalty0 (1):\penalty0 87--106, 2007.

\bibitem[Christiansen(2015)]{Christiansen:15}
S.~H. Christiansen.
\newblock Finite element systems of differential forms.
\newblock Technical Report \texttt{http://arxiv.org/abs/1006.4779v3}, arXiv,
  2015.

\bibitem[Christiansen and Winther(2008)]{Christiansen_Winther_mathcomp_2008}
S.~H. Christiansen and R.~Winther.
\newblock Smoothed projections in finite element exterior calculus.
\newblock \emph{Math. Comp.}, 77\penalty0 (262):\penalty0 813--829, 2008.

\bibitem[Ciarlet(2012)]{Ciarlet:12}
P.~Ciarlet, Jr.
\newblock {$T$}-coercivity: application to the discretization of
  {H}elmholtz-like problems.
\newblock \emph{Comput. Math. Appl.}, 64\penalty0 (1):\penalty0 22--34, 2012.

\bibitem[Ciarlet(2013)]{Ciarlet:13}
P.~Ciarlet, Jr.
\newblock Analysis of the {S}cott-{Z}hang interpolation in the fractional order
  {S}obolev spaces.
\newblock \emph{J. Numer. Math.}, 21\penalty0 (3):\penalty0 173--180, 2013.

\bibitem[Ciarlet(2016)]{Ciarlet:16}
P.~Ciarlet, Jr.
\newblock On the approximation of electromagnetic fields by edge finite
  elements. {P}art 1: {S}harp interpolation results for low-regularity fields.
\newblock \emph{Comput. Math. Appl.}, 71\penalty0 (1):\penalty0 85--104, 2016.

\bibitem[Costabel(1990)]{Costabel_1990}
M.~Costabel.
\newblock A remark on the regularity of solutions of {M}axwell's equations on
  {L}ipschitz domains.
\newblock \emph{Math. Methods Appl. Sci.}, 12\penalty0 (4):\penalty0 365--368,
  1990.

\bibitem[Ern and Guermond(2016{\natexlab{a}})]{ErnGuermond:15}
A.~Ern and J.-L. Guermond.
\newblock Finite element quasi-interpolation and best approximation.
\newblock \emph{M2AN Math. Model. Numer. Anal.}, 2016{\natexlab{a}}.
\newblock \doi{https://doi.org/10.1051/m2an/2016066}.

\bibitem[Ern and Guermond(2016{\natexlab{b}})]{Ern_Guermond_CMAM_2016}
A.~Ern and J.-L. Guermond.
\newblock Mollification in strongly {L}ipschitz domains with application to
  continuous and discrete de {R}ham complexes.
\newblock \emph{Comput. Methods Appl. Math.}, 16\penalty0 (1):\penalty0 51--75,
  2016{\natexlab{b}}.

\bibitem[Falk and Winther(2014)]{FalWi:12}
R.~S. Falk and R.~Winther.
\newblock Local bounded cochain projections.
\newblock \emph{Math. Comp.}, 83\penalty0 (290):\penalty0 2631--2656, 2014.

\bibitem[Hiptmair(1999)]{Hiptmair:99}
R.~Hiptmair.
\newblock Canonical construction of finite elements.
\newblock \emph{Math. Comp.}, 68\penalty0 (228):\penalty0 1325--1346, 1999.

\bibitem[Hiptmair(2002)]{Hiptmair_acta_numer_2002}
R.~Hiptmair.
\newblock Finite elements in computational electromagnetism.
\newblock \emph{Acta Numerica}, 11:\penalty0 237--339, 1 2002.

\bibitem[Jochmann(1999)]{Jochmann_maxwell_1999}
F.~Jochmann.
\newblock Regularity of weak solutions of {M}axwell's equations with mixed
  boundary-conditions.
\newblock \emph{Math. Methods Appl. Sci.}, 22\penalty0 (14):\penalty0
  1255--1274, 1999.

\bibitem[Kikuchi(1989)]{Kikuchi_1989}
F.~Kikuchi.
\newblock On a discrete compactness property for the {N}\'ed\'elec finite
  elements.
\newblock \emph{J. Fac. Sci. Univ. Tokyo Sect. IA Math.}, 36\penalty0
  (3):\penalty0 479--490, 1989.

\bibitem[Monk(1992)]{Monk_1992}
P.~Monk.
\newblock A finite element method for approximating the time-harmonic {M}axwell
  equations.
\newblock \emph{Numer. Math.}, 63\penalty0 (2):\penalty0 243--261, 1992.

\bibitem[Monk(2003)]{Monk_2003}
P.~Monk.
\newblock \emph{Finite element methods for {M}axwell's equations}.
\newblock Numerical Mathematics and Scientific Computation. Oxford University
  Press, New York, 2003.

\bibitem[Monk and Demkowicz(2001)]{Monk_Dem_2001}
P.~Monk and L.~Demkowicz.
\newblock Discrete compactness and the approximation of {M}axwell's equations
  in {${\mathbb R}\sp 3$}.
\newblock \emph{Math. Comp.}, 70\penalty0 (234):\penalty0 507--523, 2001.

\bibitem[Sayas(2004)]{Sayas_2004}
F.-J. Sayas.
\newblock Aubin-{N}itsche estimates are equivalent to compact embeddings.
\newblock \emph{BIT}, 44\penalty0 (2):\penalty0 287--290, 2004.

\bibitem[Sch{\"o}berl(2001)]{Schoberl_2001}
J.~Sch{\"o}berl.
\newblock Commuting quasi-interpolation operators for mixed finite elements.
\newblock Technical Report ISC-01-10-MATH, Texas A\&M University, 2001.
\newblock URL \url{www.isc.tamu.edu/publications-reports/tr/0110.pdf}.

\bibitem[Sch{\"o}berl(2005)]{Schoberl:05}
J.~Sch{\"o}berl.
\newblock A multilevel decomposition result in {$H(\mathrm{curl})$}.
\newblock In P.~Wesseling, C.~W. Oosterlee, and P.~Hemker, editors,
  \emph{Multigrid, Multilevel and Multiscale Methods, EMG 2005}, 2005.

\bibitem[Sch{\"o}berl(2008)]{Schoberl_2008}
J.~Sch{\"o}berl.
\newblock A posteriori error estimates for {M}axwell equations.
\newblock \emph{Math. Comp.}, 77\penalty0 (262):\penalty0 633--649, 2008.

\bibitem[Weber(1980)]{Weber:80}
C.~Weber.
\newblock A local compactness theorem for {M}axwell's equations.
\newblock \emph{Math. Methods Appl. Sci.}, 2\penalty0 (1):\penalty0 12--25,
  1980.

\bibitem[Xu and Zikatanov(2003)]{XuZik:03}
J.~Xu and L.~Zikatanov.
\newblock Some observations on {B}abu\v ska and {B}rezzi theories.
\newblock \emph{Numer. Math.}, 94\penalty0 (1):\penalty0 195--202, 2003.

\bibitem[Zhong et~al.(2009)Zhong, Shu, Wittum, and
  Xu]{Zhong_Shi_Wittum_Xu_2009}
L.~Zhong, S.~Shu, G.~Wittum, and J.~Xu.
\newblock Optimal error estimates for {N}edelec edge elements for time-harmonic
  {M}axwell's equations.
\newblock \emph{J. Comput. Math.}, 27\penalty0 (5):\penalty0 563--572, 2009.

\end{thebibliography}
\end{document}